\documentclass[11pt]{article}
\usepackage{algorithm}
\usepackage{algpseudocode}
\usepackage{amsfonts}
\usepackage{amsmath}
\usepackage{amssymb} % 実数Rの入力する
\usepackage{amsthm} % 定理を記述する
\usepackage[title]{appendix} % 単語"appendix"を挿入する（https://tex.stackexchange.com/questions/226481/appendix-section-title）
\usepackage{autobreak}
\usepackage{bigdelim} % 表の周りに波括弧とコメントをつける（https://motchy99.blog.fc2.com/blog-entry-78.html）
\usepackage{blkarray}
\makeatletter
\def\BA@fnsymbol#1{\ensuremath{
  \ifcase#1\or *\or \dagger\or \ddagger\or \mathsection\or \mathparagraph\or \|\or **\or \dagger\dagger\or \ddagger\ddagger \else@ctrerr\fi}}%
\makeatother
\usepackage{bm}
\usepackage{booktabs}
\usepackage{cancel} % 数式の上からキャンセル線を引く
\usepackage{cases}
\usepackage{chngcntr}
\usepackage{colortbl} % 表に色付けする
\usepackage{comment}
\usepackage{empheq} % 連立方程式を書く（https://tex.stackexchange.com/questions/451687/main-label-in-a-system-of-equations）
\usepackage{enumitem}
\setlist{leftmargin=*} % enumitemの一括設定
\usepackage{fancyhdr} % 上側に線を引く
\usepackage{fancyvrb} % verbatimを使う
\usepackage{framed} % 枠囲いをする（その２）
\usepackage{graphicx}% 図表を挿入する [dvipdfmx]
\usepackage{hyperref} % dvipdfmxが無いとリンクを貼れない
\hypersetup{
  colorlinks,
  linkcolor={red!50!black},
  citecolor={blue!50!black},
  urlcolor={blue!80!black}
}
\usepackage[capitalize]{cleveref} % hyperrefの後に読み込む
\usepackage{mathrsfs} % 花文字
\usepackage{multirow} % 表で縦方向に結合する
\usepackage[sort&compress,numbers]{natbib}
\usepackage{pifont}% http://ctan.org/pkg/pifont cmarkとxmarkで使用
\usepackage{setspace} % ダブルスペースを使う
\usepackage{stmaryrd}
\usepackage{subfig} % 図と表を並置する
\usepackage{pdflscape} % 表の横向き配置をし、用紙も横向きにする（横向きにしたい箇所を\begin/end{landscape}で挟む）（https://oxon.hatenablog.com/entry/20120428/1335624229）
\usepackage{pgfplots} % hyperrefを使用するために必要
\usepackage{url}
\usepackage{tabularx}
\usepackage[colorinlistoftodos,prependcaption,textsize=footnotesize,color=gray!40]{todonotes}

\usepackage[top=35truemm,bottom=35truemm,left=30truemm,right=30truemm]{geometry} % 余白の設定
\makeatletter

\@addtoreset{equation}{section}
\makeatother
\algtext*{EndWhile}
\algtext*{EndIf}
\algtext*{EndFor}
\theoremstyle{plain}
\newtheorem{theorem}{Theorem}[section]
\newtheorem{corollary}[theorem]{Corollary}
\newtheorem{lemma}[theorem]{Lemma}

\theoremstyle{definition}
\newtheorem{definition}[theorem]{Definition}%
\newtheorem{assumption}[theorem]{Assumption}
\theoremstyle{remark}
\newtheorem{example}[theorem]{Example}%
\DeclareMathOperator*{\cl}{cl}

\DeclareMathOperator*{\setint}{int}

\DeclareMathOperator*{\epi}{epi}

\DeclareMathOperator*{\svec}{svec}

\DeclareMathOperator*{\smat}{smat}

\DeclareMathOperator{\Sp}{Sp}

\newcommand{\relmiddle}[1]{\mathrel{}\middle#1\mathrel{}}
\newcommand{\overbar}[1]{\mkern 1.5mu\overline{\mkern-1.5mu#1\mkern-1.5mu}\mkern 1.5mu}
\newcommand{\abs}[1]{\lvert #1 \rvert} % 絶対値
\newcommand{\norm}[1]{\lVert #1 \rVert} % ノルム

\newcommand{\bbR}{\mathbb{R}}

\newcommand{\calO}{\mathcal{O}}

\newcommand{\calS}{\mathcal{S}}

\newcommand{\rmO}{\mathrm{O}}

\newcommand{\rmd}{\mathrm{d}}
\newcommand{\rme}{\mathrm{e}}

\newcommand{\rmo}{\mathrm{o}}

\newcommand{\Kexp}{K_{\rm exp}}

\makeatletter
\newcommand{\dotpplus}{\mathbin{\text{\@dotpplus}}}
\newcommand{\@dotpplus}{%
  \ooalign{\hidewidth\hbox{$\circ$}\hidewidth\cr$\m@th+$\cr}%
}
\makeatother

\makeatletter
\newcommand{\dotminus}{\mathbin{\text{\@dotminus}}}
\newcommand{\@dotminus}{%
  \ooalign{\hidewidth\hbox{$\circ$}\hidewidth\cr$\m@th-$\cr}%
}
\makeatother

\newcommand{\RNum}[1]{\uppercase\expandafter{\romannumeral #1\relax}} % ローマ数字（大文字）
\newcommand{\Rnum}[1]{\lowercase\expandafter{\romannumeral #1\relax}} % ローマ数字（小文字）

\title{Uniform symmetric spectral finiteness\\
of cones definable in o-minimal structures}

\makeatletter
\let\@fnsymbol\@arabic
\makeatother

\author{
\normalsize
    Mitsuhiro Nishijima\thanks{Department of Industrial and Systems Engineering, Keio University, 3-14-1 Hiyoshi, Kohoku-ku, Yokohama-shi, 2238522, Kanagawa, Japan ({\tt nishijima@keio.jp}).}
}

\begin{document}
\maketitle

\begin{abstract}\noindent
Cone-constrained eigenvalues of a matrix are defined by a complementarity problem over a closed convex cone.
We show that every closed convex cone that is definable in an o-minimal structure is \emph{uniformly symmetrically spectrally finite}; that is, the cone spectra of symmetric matrices relative to the cone (the sets of their cone-constrained eigenvalues) have uniformly bounded cardinality.
In particular, we derive explicit uniform bounds for semialgebraic cones, and the explicit uniform bound we obtain for a L\"{o}wner cone depends only on the order of the matrices in the cone.
Using the uniform finiteness of cone spectra, we also show that the number of singular values relative to two closed convex cones definable in a common o-minimal structure is uniformly bounded.
Finally, we provide a closed convex cone that is not definable in any o-minimal structure and for which the spectra of symmetric matrices are finite but have no uniform cardinality bound.
\end{abstract}
\vspace{0.5cm}

\vspace{0.5cm}
\noindent
{\bf MSC codes. }03C64, 14P10, 15A18, 15B48, 47A75, 90C33
% 03C64: 順序構造のモデル理論; O-ミニマリティ (Model theory of ordered structures; o-minimality)
% 14P10: 半代数的集合と関連空間 (Semialgebraic sets and related spaces)
% 15A18: 固有値、特異値、固有ベクトル (Eigenvalues, singular values, and eigenvectors)
% 15B48: 正の行列とその一般化; 行列の円錐 (Positive matrices and their generalizations; cones of matrices)
% 47A75: 線形演算子の固有値問題 (Eigenvalue problems for linear operators)
% 90C20: Quadratic programming
% 90C22: Semidefinite programming
% 90C25: Convex programming
% 90C26: Nonconvex programming, global optimization
% 90C33: 相補性および平衡問題、変分不等式 (有限次元) (数理計画法の側面) (Complementarity and equilibrium problems and variational inequalities (finite dimensions) (aspects of mathematical programming))

\section{Introduction}\label{sec:intro}
Let $K$ be a closed convex cone in $n$-dimensional Euclidean space and let $\bm{A}$ be a real $n\times n$ matrix.
A real number $\lambda$ is called a \emph{$K$-eigenvalue} of $\bm{A}$ if there exists a nonzero $\bm{x} \in K$ such that $\bm{A}\bm{x} - \lambda\bm{x}$ belongs to the dual cone of $K$ and is orthogonal to $\bm{x}$.
Cone-constrained eigenvalues arise frequently in applied mathematics, especially in mathematical optimization.
For example, they are closely tied to \emph{copositivity over $K$} (also referred to, e.g., as \emph{$K$-copositivity})~\cite{HS2010}.
If $K$ is nonzero and $\bm{A}$ is symmetric, then $\bm{A}$ is copositive over $K$ if and only if its minimum $K$-eigenvalue is nonnegative; see \cite[Equation~(8)]{ST2009} and \cite[Corollary~6.2]{Gow2012}.
Moreover, cone-constrained eigenvalues play a role in the study of Rayleigh quotient minimization under a conic constraint~\cite{ST2009} and in spectral graph theory~\cite{FJT2017}.

In this paper, we address the finiteness of the \emph{$K$-spectra} of symmetric matrices, that is, the sets of their $K$-eigenvalues.
As discussed in \cite[Question~2]{ST2003}, \cite[page~773]{ZG2009}, and \cite{See2020}, a central problem is to identify the class of closed convex cones for which the associated cone spectrum is finite for every symmetric matrix.
A closed convex cone satisfying this finiteness property is called \emph{symmetrically spectrally finite} in \cite{See2020}.
Examples of symmetrically spectrally finite cones are summarized in \cite[Table~1]{See2020}.
Specifically, polyhedral cones~\cite[Proposition~6.1]{See1999}, Lorentz cones (also called second-order cones)~\cite[Corollary~4.5]{ST2003} (see also \cite[Proposition~3]{See2020}), and certain cones that are linearly isomorphic to a Lorentz cone~\cite[Propositions~5 and 6]{See2020} are symmetrically spectrally finite.

It was unknown whether L\"{o}wner cones (cones of positive semidefinite matrices) and direct products of Lorentz cones are symmetrically spectrally finite, as noted in \cite[Table~1]{See2020}.
These cones belong to the class of \emph{symmetric cones}~\cite{FK1994}, that is, self-dual and homogeneous cones.
The case of symmetric cones is of particular interest in the context of copositivity, as copositivity over these cones has been investigated intensively~\cite{NN2024_App,NN2024_Gen,NL2025,Nis2026,NL2026,FGN+2024,HX20XX,Orl2021}.
Although previous studies left open whether all symmetric cones are symmetrically spectrally finite, the finiteness of the associated cone spectra of certain matrices was shown in \cite[Theorem~12]{ZG2009} and \cite[Theorem~16]{GTM2009}.
 
It should be noted that symmetric spectral finiteness concerns only symmetric matrices.
The stronger property that the cone spectrum is finite for every (not necessarily symmetric) matrix occurs less frequently, although examples are known, as shown in \cite[Proposition~6.1]{See1999} and \cite[Proposition~6]{See2020}.
Indeed, a spectrum relative to the $3$-dimensional Lorentz cone may contain a nondegenerate interval~\cite[Example~3]{Qui1986}.

In this paper, we prove that closed convex cones that are \emph{definable in an o-minimal structure} are symmetrically spectrally finite.
In fact, such cones are shown to be \emph{uniformly} symmetrically spectrally finite (\cref{thm:omin_ussf}) in that the cardinality of the associated cone spectrum of each symmetric matrix is bounded above by a constant that is independent of the matrix.

Roughly speaking, a set is definable if it can be described in a tame and well-behaved manner; the precise definition is provided in \cref{subsec:omin}.
Typical examples of definable sets are semialgebraic sets, i.e., sets that can be expressed as a finite union of sets, each defined by finitely many polynomial equalities and inequalities.
For example, polyhedral cones and Lorentz cones are semialgebraic.
Since they are uniformly symmetrically spectrally finite by \cite[Proposition~6.1]{See1999} and \cite[Proposition~3]{See2020}, respectively, our theorem recovers the known uniform finiteness property in both cases.
Other examples of semialgebraic cones listed in \cite[Table~1]{See2020} include symmetric cones~\cite[Section~\RNum{3}.2]{FK1994} and cones of completely positive matrices~\cite{PS2024}.
Thus, they are in fact uniformly symmetrically spectrally finite.
Cones involving the exponential or logarithmic function, such as the exponential cone, can also be definable. 

We obtain two further results as consequences of \cref{thm:omin_ussf}.
First, in \cref{cor:singlar_val}, we show that the number of singular values relative to two closed convex cones definable in a common o-minimal structure is uniformly bounded.
Cone-constrained singular values were introduced in \cite{SS2023_Con}, and the finiteness of their number is also discussed in \cite{SS2023_Con,SS2023_Sin}.
Second, in \cref{thm:semialg_uniform_bounded}, we provide an explicit uniform bound for the cardinalities of spectra of symmetric matrices relative to a closed semialgebraic convex cone when polynomial representations of the cone and its dual are given.
In particular, in \cref{ex:psd}, we provide an explicit uniform bound for the cardinalities of spectra of symmetric matrices relative to a (vectorized) L\"{o}wner cone.
This bound depends only on the order of the positive semidefinite matrices defining the cone.

A natural question that arises from \cref{thm:omin_ussf} concerns the (uniform) symmetric spectral finiteness of cones that are not definable in any o-minimal structure. 
Undefinable cones are not guaranteed to be symmetrically spectrally finite, and (necessarily undefinable) cones that are not symmetrically spectrally finite are given in \cite[Section~6]{IS2007} and \cite[Theorem~6]{See2021}.
In \cref{subsec:nonsemialg_finite_nonunif}, we provide an undefinable closed convex cone $K$ that is symmetrically spectrally finite but is not uniformly symmetrically spectrally finite.
In other words, the $K$-spectrum is finite for every symmetric matrix, but for any prescribed bound we can choose a symmetric matrix whose $K$-spectrum has cardinality at least that bound.

This paper is organized as follows.
In \cref{sec:preliminaries}, we introduce the notation and results used in the subsequent discussion.
In \cref{sec:semialg_finite}, we show the main result, i.e., the uniform symmetric spectral finiteness of definable cones.
In \cref{sec:appl}, we present two applications of the main result.
In \cref{sec:nonsemialg}, we illustrate the uniform symmetric spectral finiteness of the exponential cone and provide an undefinable closed convex cone that is symmetrically spectrally finite but is not uniformly symmetrically spectrally finite.

\section{Preliminaries}\label{sec:preliminaries}
\subsection{Basic notation}
$\bbR^{m\times n}$ denotes the space of real $m\times n$ matrices and $\calS^n$ denotes the space of symmetric matrices in $\bbR^{n\times n}$.
Boldface uppercase letters, such as $\bm{A}$, denote matrices in $\bbR^{m\times n}$.
For a matrix $\bm{A}$, the transpose of $\bm{A}$ is denoted by $\bm{A}^\top$ and the $(i,j)$th entry of $\bm{A}$ is denoted by $A_{ij}$.
We write $\det\bm{A}$ for the determinant of a square matrix $\bm{A}$.
We use $\bm{I}$ to denote the identity matrix, whose size is determined by the context, and write $\bbR\bm{I}$ for the set of real scalar multiples of $\bm{I}$.
We write $\bbR^{n\times 1}$ as $\bbR^n$ and equip it with the standard inner product and the induced norm denoted by $\norm{\cdot}$.
We use $\bbR_+^n$ (resp. $\bbR_{++}^n$) to denote the set of entrywise nonnegative (resp. positive) vectors in $\bbR^n$.
Boldface lowercase letters, such as $\bm{a}$, denote vectors in $\bbR^n$, and $a_i$ denotes the $i$th entry of $\bm{a}$.
We use $\bm{e}_i$ to denote the vector whose $i$th entry is $1$ and whose other entries are $0$.
For $\bm{a},\bm{b}\in \bbR^n$, we use $\bm{a} \perp \bm{b}$ to represent $\bm{a}^\top\bm{b} = 0$.

For a real number $a$, we use $\abs{a}$ to denote the absolute value of $a$.
All polynomials are assumed to have real coefficients in this paper, and the degree of a polynomial $f$ is denoted by $\deg f$.
For two functions $\bm{f}$ and $\bm{g}$, we use $\bm{f}\circ \bm{g}$ to denote the composite function.
For a univariate scalar-valued function $f$, we write $f'$ for the derivative of $f$.
For a univariate vector-valued function $\bm{f} = (f_1,\dots,f_n)^\top$, we define $\bm{f}' \coloneqq (f_1',\dots,f_n')^\top$.
For a function $f\colon \bbR^n \to (-\infty,+\infty]$, we define
\begin{equation*}
\epi f \coloneqq \{(\bm{x}^\top,y)^\top \in \bbR^{n+1} \mid f(\bm{x})\le y\}.
\end{equation*}

For a positive integer $n$, we define $[n] \coloneqq \{1,\dots,n\}$.
For a set $S$, we use $\#S$ to denote the number of elements in $S$.
A set $S$ in $\bbR^n$ is called a \emph{cone} if $a\bm{x} \in S$ for all $a \ge 0$ and $\bm{x} \in S$.
For a set $S$ in $\bbR^n$, we use $\setint S$, $\cl S$, and $S^*$ to denote the interior of $S$, the closure of $S$, and the dual cone of $S$, respectively.
The dual cone is the cone of $\bm{x}\in \bbR^n$ such that $\bm{x}^\top\bm{y}\ge 0$ for all $\bm{y}\in S$.
For a topological space $S$, the number of connected components of $S$ is denoted by $b_0(S)$.

\subsection{O-minimal structures and definable sets}\label{subsec:omin}
This subsection provides preliminaries on o-minimal structures on the field $\bbR$.
O-minimal structures are defined in several equivalent ways; see, for example, \cite[Definition~6]{BDL+2007}, \cite[Definition~1.4]{Cos1999}, and \cite[page~2]{van1998}.
We follow \cite[Definition~6]{BDL+2007} in giving the definition below.

\begin{definition}\label{def:omin}
An \emph{o-minimal structure} (on the field $\bbR$) is a sequence $(\calO_n)_{n\ge 1}$ satisfying the following properties for each $n\ge 1$:
\begin{enumerate}[label=(\roman*), ref=\roman*]
\item $\calO_n$ is a Boolean subalgebra of the power set of $\bbR^n$;
\item if $S \in \calO_n$, then $\bbR\times S \in \calO_{n+1}$ and $S \times \bbR \in \calO_{n+1}$;
\item if $S \in \calO_{n+1}$, then the projection of $S$ onto the first $n$ coordinates belongs to $\calO_n$; \label{enum:omin_proj}
\item $\calO_n$ contains all sets of the form $\{\bm{x}\in \bbR^n \mid f(\bm{x}) = 0\}$, where $f$ is a polynomial;
\item $\calO_1$ is precisely the family of finite unions of intervals and points.
\end{enumerate}
\end{definition}

Let $\mathcal{O}$ be a fixed o-minimal structure.
Every element in $\calO_n$ with $n\ge 1$ is called \emph{definable} (in $\calO$).
A function $\bm{f}\colon S \to \bbR^m$ is called \emph{definable} if the graph of $\bm{f}$ is definable.

The semialgebraic sets form a standard example of an o-minimal structure.
A set $S$ in $\bbR^n$ is called \emph{semialgebraic} if there exist polynomials $f_{ik}$ in $n$ variables and relations $*_{ik} \in \{<,=\}$ for $i\in [N]$ and $k\in [s_i]$ such that
\begin{equation*}
S = \bigcup_{i=1}^N \{\bm{x}\in \bbR^n \mid f_{ik}(\bm{x}) *_{ik} 0 \text{ for all $k\in [s_i]$}\}.
\end{equation*}
As a consequence of the Tarski--Seidenberg theorem~\cite[Theorem~2.2.1]{BCR1998}, the collection of semialgebraic sets constitutes an o-minimal structure.

Some o-minimal structures also contain non-semialgebraic sets defined using the exponential function.
More precisely, for each $n\ge 1$, we define $\calO_n$ as the family of sets of the form
\begin{equation*}
\pi(\{(\bm{x},\bm{y}) \in \bbR^n \times \bbR^k \mid f(x_1,\dots,x_n,y_1,\dots,y_k,\rme^{x_1},\dots,\rme^{x_n},\rme^{y_1},\dots,\rme^{y_k}) = 0\}),
\end{equation*}
where $k\ge 0$, $\pi\colon \bbR^n \times \bbR^k \to \bbR^n$ is the projection, and $f$ is a polynomial in $2(n+k)$ variables.
Then $(\calO_n)_{n\ge 1}$ forms an o-minimal structure; see \cite{Wil1996} and \cite[page~3]{van1998}.
Cones that are non-semialgebraic but definable in this structure also arise in optimization.
For example, the exponential cone~\cite[Section~4.2]{Cha2009}
\begin{equation}
\Kexp \coloneqq \{\bm{x}\in \bbR^3 \mid x_2 > 0,\ x_2\rme^{x_1/x_2} \le x_3 \} \cup \{-a\bm{e}_1 + b\bm{e}_3 \mid a,b\ge 0\}, \label{eq:def_Kexp}
\end{equation}
which we will revisit in \cref{subsec:expcone}, is definable.
Furthermore, using the definability of the logarithmic function (the inverse of the exponential function), together with the identity $x^a=\rme^{a\log(x)}$ for $a\in\bbR$ and $x>0$, we see that definable closed convex cones include
\begin{enumerate}[label=(\roman*), ref=\roman*]
\item the $p$-order cone $\{\bm{x}\in \bbR^n \mid x_1 \ge (\abs{x_2}^p + \cdots + \abs{x_n}^p)^{1/p}\}$, where $p \ge 1$ and $n\ge 2$~\cite[Section~2.2.3]{BV2004},
\item the log-determinant cone $\cl\{(x,y,\bm{Z}) \in \bbR\times \bbR_{++}\times \calS_{++}^d \mid x \le y\log\det(\bm{Z}/y)\}$, where $d \ge 1$ and $\calS_{++}^d$ denotes the set of positive definite matrices in $\calS^d$~\cite[Section~5.6]{CKV2022},
\item the generalized power cone $\{(\bm{x},\bm{y}) \in \bbR^m \times \bbR_+^n \mid \norm{\bm{x}} \le y_1^{\alpha_1}\cdots y_n^{\alpha_n}\}$, where $m\ge 1$, $n\ge 1$, and the weights $\alpha_1,\dots,\alpha_n > 0$ sum to $1$~\cite[Section~3.1.2]{Cha2009}, and
\item the relative entropy cone $\cl\{(\bm{x},\bm{y},z) \in \bbR_{++}^d \times \bbR_{++}^d \times \bbR \mid \sum_{i=1}^d x_i\log(x_i/y_i) \le z\}$, where $d\ge 1$~\cite[page~75]{CKV2023}.
\end{enumerate}

We summarize some known results about definable sets and functions needed for the subsequent discussion.
In the following theorems, we say that a set is \emph{definably connected} if it is definable and cannot be written as the union of two disjoint nonempty definable subsets that are open in the set.

\begin{theorem}\label{thm:omin_cc}
Let $S$ be a definable set in $\bbR^m \times \bbR^n$.
Then there exists $M > 0$ such that for every $\bm{a}\in \bbR^m$, the definable set $S(\bm{a}) \coloneqq \{\bm{x}\in \bbR^n \mid (\bm{a},\bm{x})\in S\}$ has at most $M$ connected components, each of which is definable.
\end{theorem}

\begin{proof}
For every $\bm{a}\in\bbR^m$, each connected component of $S(\bm{a})$ is definable by \cite[page~4]{van1998}.
Moreover, a definable set is connected if and only if it is definably connected~\cite[Exercise~2.19.7 on page~59]{van1998}.
Consequently, the connected components of $S(\bm{a})$ are precisely its maximal definably connected subsets.
The number of these subsets is bounded uniformly in $\bm{a}$ by \cite[Corollary~3.6 on page~60]{van1998}, and this proves the assertion.
\end{proof}

\begin{theorem}\label{thm:omin_connected}
Let $S$ be a connected definable subset of $\bbR^n$.
Then for any $\bm{x}_0,\bm{x}_1 \in S$, there exists a definable continuous curve $\bm{\gamma}:[0,1]\to S$ such that $\bm{\gamma}(0) = \bm{x}_0$ and $\bm{\gamma}(1) = \bm{x}_1$. 
\end{theorem}

\begin{proof}
Note again that, by \cite[Exercise~2.19.7 on page~59]{van1998}, $S$ is definably connected.
Therefore, the conclusion follows from \cite[Proposition~3.2 on page~100]{van1998}.
\end{proof}

\begin{theorem}\label{thm:omin_C1}
Let $I$ be an open interval of $\bbR$ and let $\bm{f}:I\to \bbR^n$ be a definable function.
Then there exists a finite subset $T$ of $I$ such that $\bm{f}$ is of class $C^1$ on $I\setminus T$.
\end{theorem}

\begin{proof}
The statement for $n = 1$ follows from \cite[Proposition~2.5 on page~109]{van1998}.
(The statement in \cite{van1998} mentions only differentiability, but its proof implies $C^1$-smoothness; see also \cite[Corollary~6.5]{Cos1999}.)
We consider the case where $n\ge 2$.
For every $i\in [n]$, since $f_i\colon I\to \bbR$ is definable by \cite[Exercise~2.5.1 on page~15]{van1998}, it follows from the result for $n = 1$ that there exists a finite subset $T_i$ of $I$ such that $f_i$ is of class $C^1$ on $I \setminus T_i$.
Letting $T \coloneqq \bigcup_{i=1}^n T_i$, we see that $T$ is also a finite subset of $I$ and $f_1,\dots,f_n$ are of class $C^1$ on $I \setminus T$.
Therefore, we obtain the desired result.  
\end{proof}

\subsection{Cone-constrained eigenvalues}
Let $K$ be a closed convex cone in $\bbR^n$ and let $\bm{A} \in \bbR^{n\times n}$.
We call $\lambda \in \bbR$ a \emph{$K$-eigenvalue of $\bm{A}$} if there exists a nonzero $\bm{x}$ such that $K \ni \bm{x} \perp (\bm{A}\bm{x} - \lambda\bm{x})\in K^*$.
Such a vector $\bm{x}$ is referred to as a \emph{$K$-eigenvector of $\bm{A}$} corresponding to $\lambda$.
The set of $K$-eigenvalues of $\bm{A}$ is denoted by $\Sp(\bm{A},K)$ and called the \emph{$K$-spectrum of $\bm{A}$} or the \emph{spectrum of $\bm{A}$ relative to $K$}.

We recall the notion of finiteness of $K$-spectra of symmetric matrices and introduce its uniform variant.
\begin{definition}
Let $K$ be a closed convex cone  in $\bbR^n$.
\begin{enumerate}[label=(\roman*), ref=\roman*]
\item $K$ is called \emph{symmetrically spectrally finite} if $\Sp(\bm{A},K)$ is finite for any $\bm{A}\in \calS^n$~\cite[Definition~1.d]{See2020}.
\item $K$ is called \emph{uniformly symmetrically spectrally finite} if there exists $M > 0$ such that $\#\Sp(\bm{A},K) \le M$ for any $\bm{A}\in \calS^n$.
\end{enumerate}
\end{definition}
As can be seen from the definitions, symmetric spectral finiteness is implied by uniform symmetric spectral finiteness.
The existence of a closed convex cone that is symmetrically spectrally finite but is not uniformly symmetrically spectrally finite is nontrivial.
We construct such a cone in \cref{subsec:nonsemialg_finite_nonunif}.

The following lemma describes the $K$-spectrum of a multiple of the identity matrix.
The proof is straightforward, so we omit it.

\begin{lemma}\label{lem:property_spec}
Let $K$ be a nonzero closed convex cone in $\bbR^n$ and let $\lambda\in \bbR$.
Then we have $\Sp(\lambda\bm{I},K) = \{\lambda\}$.
\end{lemma}

\section{Uniform symmetric spectral finiteness of definable cones}\label{sec:semialg_finite}
In this section, we prove our main result, namely that definable closed convex cones are uniformly symmetrically spectrally finite.
We introduce some additional notation for the following theorem.
For $\bm{A} \in \calS^n$ and $\bm{x}\in \bbR^n$, we define
\begin{align*}
q_{\bm{A}}(\bm{x}) &\coloneqq \bm{x}^\top\bm{A}\bm{x}, \nonumber\\
\bm{s}_{\bm{A}}(\bm{x}) &\coloneqq \bm{A}\bm{x} - q_{\bm{A}}(\bm{x})\bm{x}.
\end{align*}
In addition, for a closed convex cone $K$ in $\bbR^n$, we define
\begin{equation*}
X(\bm{A}) \coloneqq \{\bm{x} \in K \mid \norm{\bm{x}} = 1,\ \bm{s}_{\bm{A}}(\bm{x}) \in K^*\} = \left\{\bm{x} \in K \relmiddle| \sum_{i=1}^nx_i^2 =1,\ \bm{s}_{\bm{A}}(\bm{x}) \in K^*\right\}.
\end{equation*}
This is the set of normalized $K$-eigenvectors of $\bm{A}$, and thus we have
\begin{equation}
\Sp(\bm{A},K) = q_{\bm{A}}(X(\bm{A})). \label{eq:qAXA}
\end{equation}

\begin{theorem}\label{thm:omin_ussf}
Let $K$ be a closed convex cone in $\bbR^n$ that is definable in an o-minimal structure.
Then there exists $M > 0$ such that for every $\bm{A} \in \calS^n$, we have
\begin{equation*}
\#\Sp(\bm{A},K) \le b_0(X(\bm{A})) \le M.
\end{equation*}
In particular, $K$ is uniformly symmetrically spectrally finite.
\end{theorem}

\begin{proof}
Since $K$ is definable, it follows from \eqref{enum:omin_proj} of \cref{def:omin} (see also \cite[Theorem~1.13]{Cos1999}) that $K^*$ is definable and so is the set
\begin{equation*}
X \coloneqq \{(\bm{A},\bm{x}) \in \calS^n \times K \mid \norm{\bm{x}} = 1,\ \bm{s}_{\bm{A}}(\bm{x}) \in K^*\}.
\end{equation*}
By \cref{thm:omin_cc}, there exists $M > 0$ such that for every $\bm{A} \in \calS^n$, there are at most $M$ connected components of $X(\bm{A}) = \{\bm{x} \in \bbR^n \mid (\bm{A},\bm{x}) \in X\}$, each of which is definable.
Let $k\coloneqq b_0(X(\bm{A}))$ for simplicity, and denote these components by $X_1,\dots,X_k$.

For every $i \in [k]$, the function $q_{\bm{A}}|_{X_i}$ is constant.
Indeed, since $X_i$ is connected and definable, \cref{thm:omin_connected} implies that for any $\bm{x}_0,\bm{x}_1\in X_i$, there exists a definable continuous curve $\bm{\gamma} \colon [0,1] \to X_i$ such that $\bm{\gamma}(0) = \bm{x}_0$ and $\bm{\gamma}(1) = \bm{x}_1$.
By \cref{thm:omin_C1}, there exist $t_0,\dots,t_m\in [0,1]$ with $0 = t_0 < \cdots < t_m = 1$ such that $\bm{\gamma}|_{(t_{j-1},t_j)}$ is of class $C^1$ for every $j \in [m]$.
We define $f\colon [0,1]\to \bbR$ as $f \coloneqq q_{\bm{A}} \circ \bm{\gamma}$.
We claim that, for each $j\in [m]$, the function $f|_{(t_{j-1},t_j)}$ is constant.
To prove the claim, we define the curve $\bm{c}\colon [0,1] \to K^*$ as $\bm{c} \coloneqq \bm{s}_{\bm{A}} \circ \bm{\gamma}$, where $\bm{c}(t) \in K^*$ for $t\in [0,1]$ follows from $\bm{\gamma}(t)\in X_i$, $X_i \subseteq X(\bm{A})$, and the definition of $X(\bm{A})$.
Let $\tau \in (t_{j-1},t_j)$ be arbitrary.
For any $t\in (t_{j-1},t_j)$, since $\bm{\gamma}(t) \in K$ and $\bm{c}(\tau) \in K^*$, we have $\bm{c}(\tau)^\top\bm{\gamma}(t) \ge 0$.
Moreover, since $\norm{\bm{\gamma}(\tau)} = 1$, we have $\bm{c}(\tau)^\top\bm{\gamma}(\tau) = 0$.
This means that the function $\bm{c}(\tau)^\top\bm{\gamma}(t)\colon (t_{j-1},t_j)\to \bbR$ attains its global minimum at $t = \tau$, and thus we have
\begin{equation}
0 = \left.\frac{\rmd}{\rmd t}\bm{c}(\tau)^\top\bm{\gamma}(t)\right|_{t = \tau}  = \bm{c}(\tau)^\top\bm{\gamma}'(\tau). \label{eq:ldotg_diff_tau}
\end{equation}
In addition, since $\bm{\gamma}(t)^\top\bm{\gamma}(t) = 1$ for all $t\in [0,1]$, we have
\begin{equation}
0 = \left.\frac{\rmd}{\rmd t}\bm{\gamma}(t)^\top\bm{\gamma}(t)\right|_{t = \tau}  = 2\bm{\gamma}(\tau)^\top\bm{\gamma}'(\tau). \label{eq:gdotg_diff_tau}
\end{equation}
Therefore, it follows that
\begin{equation}
f'(\tau) = 2(\bm{A}\bm{\gamma}(\tau))^\top\bm{\gamma}'(\tau) = 2(\bm{c}(\tau) + q_{\bm{A}}(\bm{\gamma}(\tau))\bm{\gamma}(\tau))^\top \bm{\gamma}'(\tau) = 0, \label{eq:fd_tau}
\end{equation}
where we use the symmetry of $\bm{A}$ to derive the first equality, use the definition of $\bm{c}$ to derive the second equality, and use \eqref{eq:ldotg_diff_tau} and \eqref{eq:gdotg_diff_tau} to derive the third equality.
Since $\tau \in (t_{j-1},t_j)$ is arbitrary, $f|_{(t_{j-1},t_j)}$ is constant.
Combining this with the continuity of $f$, we see that $f$ is constant.
In particular, it follows that $f(0) = f(1)$, i.e., $q_{\bm{A}}(\bm{x}_0) = q_{\bm{A}}(\bm{x}_1)$.
Since $\bm{x}_0,\bm{x}_1\in X_i$ are arbitrary, we conclude that $q_{\bm{A}}|_{X_i}$ is constant.

For every $i \in [k]$, we write $q_i$ for the value of $q_{\bm{A}}$ on $X_i$.
Since the union of $X_1,\dots,X_k$ is $X(\bm{A})$, it follows from \eqref{eq:qAXA} that $\Sp(\bm{A},K) = \{q_1,\dots,q_k\}$.
Therefore, we have $\#\Sp(\bm{A},K) \le k$.
Since $k \le M$ and $M$ is independent of $\bm{A}$, we obtain the desired result.
\end{proof}

We note that the symmetry of $\bm{A}$, which was used in \eqref{eq:fd_tau}, is essential for \cref{thm:omin_ussf}.
If $\bm{A}$ is not symmetric, we cannot apply \cref{thm:omin_ussf}, and the finiteness of $\Sp(\bm{A},K)$ is not guaranteed.
For example, Quittner~\cite[Example~3]{Qui1986} provided a $3\times 3$ nonsymmetric matrix whose spectrum relative to the $3$-dimensional Lorentz cone contains a nondegenerate interval; see also \cite{See2020,See2021} for other classes of cones that admit an infinite spectrum for nonsymmetric matrices.

As discussed in \cref{sec:intro,subsec:omin}, numerous closed convex cones of interest in optimization are definable and thus uniformly symmetrically spectrally finite.
Because constructing a counterexample to the converse of \cref{thm:omin_ussf} is rather technical, we do not pursue it here.

\section{Applications}\label{sec:appl}
In this section, we derive two results as consequences of \cref{thm:omin_ussf}.
In \cref{subsec:singval}, we show that the number of singular values relative to two closed convex cones definable in a common o-minimal structure is uniformly bounded.
In \cref{subsec:semialg}, we provide an explicit uniform bound for the cardinalities of spectra of symmetric matrices relative to a closed semialgebraic convex cone when polynomial representations of the cone and its dual are given.

\subsection{Uniform finiteness of singular value spectra relative to definable cones} \label{subsec:singval}
Let $P$ and $Q$ be closed convex cones in $\bbR^m$ and $\bbR^n$, respectively, and let $\bm{A}\in \bbR^{m\times n}$.
We call $\sigma \in \bbR$ a \emph{$(P,Q)$-singular value} of $\bm{A}$~\cite[Definition~1.1]{SS2023_Con} if there exist $\bm{u} \in \bbR^m$ and $\bm{v}\in \bbR^n$ with $\norm{\bm{u}} = \norm{\bm{v}} = 1$ such that $P \ni \bm{u} \perp (\bm{A}\bm{v} - \sigma\bm{u}) \in P^*$ and $Q \ni \bm{v} \perp (\bm{A}^\top\bm{u} - \sigma\bm{v}) \in Q^*$.

In the following corollary, we show that if $P$ and $Q$ are definable in a common o-minimal structure, then for any matrix $\bm{A} \in \bbR^{m\times n}$, the set of $(P,Q)$-singular values of $\bm{A}$, denoted by $\Sigma(\bm{A},P,Q)$, has cardinality bounded above by a constant that is independent of $\bm{A}$.

\begin{corollary}\label{cor:singlar_val}
Let $P$ and $Q$ be closed convex cones in $\bbR^m$ and $\bbR^n$, respectively, that are definable in a common o-minimal structure.
Then there exists $M > 0$ such that $\#\Sigma(\bm{A},P,Q) \le M$ for any $\bm{A} \in \bbR^{m\times n}$.
\end{corollary}

\begin{proof}
Since $P\times Q$ is definable, it follows from \cref{thm:omin_ussf} that $P\times Q$ is uniformly symmetrically spectrally finite; namely, there exists $M > 0$ such that
\begin{equation}
\#\Sp(\bm{B},P\times Q) \le M \label{eq:PtimesQ_USSF}
\end{equation}
for any $\bm{B} \in \calS^{m+n}$.
Now, for any $\bm{A} \in \bbR^{m\times n}$, let $\bm{E}_{\bm{A}} \coloneqq \begin{psmallmatrix}
\bm{O} & \bm{A}\\
\bm{A}^\top &  \bm{O}
\end{psmallmatrix} \in \calS^{m+n}$, where $\bm{O}$ denotes the zero matrix of the appropriate size.
Then it follows from \cite[Proposition~2.5]{SS2023_Con} that
\begin{equation}
\Sigma(\bm{A},P,Q) \subseteq \Sp(\bm{E}_{\bm{A}},P\times Q). \label{eq:inclu_sig_eig}
\end{equation}
Since $\bm{E}_{\bm{A}}$ is symmetric, it follows from \eqref{eq:PtimesQ_USSF} and \eqref{eq:inclu_sig_eig} that $\#\Sigma(\bm{A},P,Q) \le M$.
Thus, we obtain the desired result.
\end{proof}

\subsection{Explicit uniform upper bound for semialgebraic cones}\label{subsec:semialg}
\cref{thm:omin_ussf} shows that an upper bound on $b_0(X(\bm{A}))$ also bounds $\#\Sp(\bm{A},K)$, but it is difficult to obtain an explicit bound for $b_0(X(\bm{A}))$ in general o-minimal structures.
In this subsection, we consider the case where $K$ is semialgebraic to derive a uniform upper bound for $b_0(X(\bm{A}))$ that is explicit in terms of polynomial representations of $K$ and $K^*$.

Let $K$ be a closed semialgebraic cone in $\bbR^n$.
Since $K$ and $K^*$ are semialgebraic, they can each be represented as a finite union of basic closed semialgebraic sets~\cite[Theorem~2.7.2]{BCR1998}; namely,
\begin{align}
K &= \bigcup_{i=1}^{N_1} \{\bm{x} \in \bbR^n \mid f_{ik}(\bm{x}) \ge 0 \text{ for all $k\in [s_i]$}\}, \label{eq:K_semialg}\\
K^* &= \bigcup_{j=1}^{N_2} \{\bm{x} \in \bbR^n \mid g_{jl}(\bm{x}) \ge 0 \text{ for all $l\in [t_j]$}\}, \label{eq:dualK_semialg}
\end{align}
where all $f_{ik}$ and $g_{jl}$ are polynomials whose degrees are assumed to be at least $1$ without loss of generality.
Using \eqref{eq:K_semialg} and \eqref{eq:dualK_semialg}, we can represent the set $X(\bm{A})$ as
\begin{equation}
X(\bm{A}) = \bigcup_{i=1}^{N_1}\bigcup_{j=1}^{N_2} \underbrace{\left\{\bm{x}\in \bbR^n \relmiddle|
\begin{aligned}
&f_{ik}(\bm{x}) \ge 0\ \text{for all $k\in [s_i]$},\\
&(g_{jl}\circ \bm{s}_{\bm{A}})(\bm{x}) \ge 0\ \text{for all $l \in [t_j]$},\\
&\sum_{k=1}^n x_k^2 - 1 \ge 0,\ 1 - \sum_{k=1}^n x_k^2 \ge 0
\end{aligned}
\right\}}_{\eqqcolon X_{ij}}. \label{eq:selfdual_S_semialg}
\end{equation}

Numerous studies give bounds on the number of connected components of semialgebraic sets; see \cite[Remark~11.5.5]{BCR1998}, \cite[Chapter~7]{BPR2006}, \cite[Section~37.7]{BM2017}, and the references therein.
We use the result of \cite{Mil1964} to bound $b_0(X(\bm{A}))$.

\begin{lemma}\label{lem:b0_Mil1964}
Let $K$ be a closed semialgebraic convex cone in $\bbR^n$, and assume that $K$ and $K^*$ are represented as in \eqref{eq:K_semialg} and \eqref{eq:dualK_semialg}, respectively.
For every $i \in [N_1]$ and $j \in [N_2]$, we define
\begin{equation*}
\overbar{D}_{ij} \coloneqq 4 + \sum_{k=1}^{s_i}\deg f_{ik} + 3\sum_{l=1}^{t_j}\deg g_{jl}.
\end{equation*}
Then for every $\bm{A}\in \calS^n$, we have
\begin{equation*}
b_0(X(\bm{A})) \le \frac{1}{2}\sum_{i=1}^{N_1}\sum_{j=1}^{N_2}(\overbar{D}_{ij} + 2)(\overbar{D}_{ij} + 1)^{n-1}.
\end{equation*}
\end{lemma}

\begin{proof}
For every $i\in [N_1]$ and $j \in [N_2]$, we write $D_{ij}$ for the sum of the degrees of the nonzero polynomials used in the representation of $X_{ij}$ defined in \eqref{eq:selfdual_S_semialg}.
By \cite[Theorem~3]{Mil1964}, we have
\begin{equation*}
b_0(X_{ij}) \le \frac{1}{2}(D_{ij} + 2)(D_{ij} + 1)^{n-1}.
\end{equation*}
Since $b_0(X(\bm{A})) \le \sum_{i=1}^{N_1}\sum_{j=1}^{N_2} b_0(X_{ij})$ and every $D_{ij}$ is no greater than $\overbar{D}_{ij}$, we obtain the desired result.
\end{proof}

From \cref{thm:omin_ussf,lem:b0_Mil1964}, we obtain the following explicit uniform upper bound.

\begin{theorem}\label{thm:semialg_uniform_bounded}
Let $K$ be a closed semialgebraic convex cone in $\bbR^n$, and assume that $K$ and $K^*$ are represented as in \eqref{eq:K_semialg} and \eqref{eq:dualK_semialg}, respectively.
Then for every $\bm{A}\in \calS^n$, we have
\begin{equation*}
\#\Sp(\bm{A},K) \le \frac{1}{2}\sum_{i=1}^{N_1}\sum_{j=1}^{N_2}(\overbar{D}_{ij} + 2)(\overbar{D}_{ij} + 1)^{n-1},
\end{equation*}
where each $\overbar{D}_{ij}$ is defined in \cref{lem:b0_Mil1964}.
\end{theorem}

A closed convex cone $K$ in $\bbR^n$ is said to be \emph{self-dual} if $K = K^*$.
In the following corollary, we derive simpler uniform bounds for the cardinalities of spectra of symmetric matrices relative to general semialgebraic self-dual cones.

\begin{corollary}\label{cor:b0_Mil1964_selfdual}
Let $K$ be a semialgebraic self-dual cone in $\bbR^n$, and assume that it is represented as in \eqref{eq:K_semialg}.
For every $i,j \in [N_1]$, we define
\begin{equation*}
\overbar{D}_{ij}' \coloneqq 4 + \sum_{k=1}^{s_i}\deg f_{ik} + 3\sum_{l=1}^{s_j}\deg f_{jl}.
\end{equation*}
Then for every $\bm{A}\in \calS^n$, we have
\begin{equation}
\#\Sp(\bm{A},K) \le \frac{1}{2}\sum_{i,j=1}^{N_1}(\overbar{D}_{ij}' + 2)(\overbar{D}_{ij}' + 1)^{n-1}. \label{eq:semialg_sefdual_ub}
\end{equation}
\end{corollary}

\begin{proof}
It follows from $K = K^*$ that $K$ and $K^*$ admit the same polynomial representation.
Since $K$ is assumed to have the polynomial representation in \eqref{eq:K_semialg}, we can also represent $K^*$ as in \eqref{eq:K_semialg}.
Therefore, by setting $N_2$ to $N_1$, $t_j$ to $s_j$ for $j\in [N_2]$, and $g_{jl}$ to $f_{jl}$ for $j\in [N_2]$ and $l \in [t_j]$ in \eqref{eq:dualK_semialg}, we obtain the desired result from \cref{thm:semialg_uniform_bounded}.
\end{proof}

We write $u(K)$ for the bound in \eqref{eq:semialg_sefdual_ub} obtained from the chosen polynomial representation of $K$.
In the following examples, we calculate explicit values of $u(K)$ when $K$ is a Pareto cone, a Lorentz cone, or a (vectorized) L\"{o}wner cone.

\begin{example}\label{ex:nno}
Let $K$ be the Pareto cone $\bbR_+^n$, which has the following polynomial representation:
\begin{equation*}
\bbR_+^n = \{\bm{x}\in \bbR^n \mid x_k \ge 0 \text{ for all $k \in [n]$}\}.
\end{equation*}
By setting $N_1 \coloneq 1$, $s_1 \coloneqq n$, and $f_{1k}(\bm{x}) \coloneqq x_k$ for $k \in [n]$, we have $\overbar{D}_{11}' = 4n+4$, and thus \cref{cor:b0_Mil1964_selfdual} yields
\begin{equation*}
u(\bbR_+^n) =  (2n+3) (4n+5)^{n-1}.
\end{equation*}
This bound grows faster than the tight uniform upper bound $2^n - 1$ for the cardinalities of $\bbR_+^n$-spectra of symmetric matrices; see \cite[Proposition~5.3]{See1999} and \cite[Proposition~3]{PS2010}.
\end{example}

\begin{example}\label{ex:Lorentz}
Let $K$ be the Lorentz cone defined as
\begin{equation*}
L_+^n \coloneqq \left\{\bm{x}\in \bbR^n \relmiddle| x_1 \ge 0,\ x_1^2 - \sum_{i=2}^n x_i^2 \ge 0\right\}.
\end{equation*}
By setting $N_1 \coloneqq 1$, $s_1 \coloneqq 2$, $f_{11}(\bm{x}) \coloneqq x_1$, and $f_{12}(\bm{x}) \coloneqq x_1^2 - \sum_{i=2}^n x_i^2$, we have $\overbar{D}_{11}' = 16$, and thus \cref{cor:b0_Mil1964_selfdual} yields 
\begin{equation*}
u(L_+^n) = 9\cdot 17^{n-1}.
\end{equation*}
This bound grows faster than the best-known uniform upper bound $3n-2$ for the cardinalities of $L_+^n$-spectra of symmetric matrices~\cite[Proposition~3]{See2020}.
\end{example}

\begin{example}\label{ex:psd}
Let $\calS_+^n$ denote the cone of positive semidefinite matrices in $\calS^n$. 
This cone is also known as the L\"{o}wner cone.
For brevity, we define $T_n \coloneqq \frac{n(n+1)}{2}$.
Let $\svec \colon \calS^n \to \bbR^{T_n}$ be the linear mapping defined by
\begin{equation*}
\svec(\bm{X}) \coloneqq (X_{11},\sqrt{2}X_{12},X_{22},\dots,\sqrt{2}X_{1n},\dots,\sqrt{2}X_{n-1,n},X_{nn})^\top
\end{equation*}
and denote the inverse of $\svec$ by $\smat\colon \bbR^{T_n} \to \calS^n$.
Then we can represent $\svec(\calS_+^n)$ as
\begin{equation*}
\svec(\calS_+^n) = \left\{\bm{x}\in \bbR^{T_n} \relmiddle| p_k(\bm{x}) \coloneqq \sum_{\substack{I\subseteq [n],\\ \#I=n-k+1}}\det\smat(\bm{x})_I \ge 0 \text{ for all $k\in [n]$}\right\},
\end{equation*}
where $\smat(\bm{x})_{I}$ is the principal submatrix obtained by extracting the rows and columns of $\smat(\bm{x})$ indexed by the elements in $I$.
See, for example, \cite[Proposition~8.2.7]{Ber2009} and \cite[Proposition~A.1]{BPT2012} for this representation.
We note that $p_k$ is a polynomial of degree $n-k+1$ for every $k \in [n]$.
By setting $N_1 \coloneqq 1$, $s_1 \coloneqq n$, and $f_{1k} \coloneqq p_k$ for $k \in [n]$, we have $\overbar{D}_{11}' = 2n^2 + 2n + 4$, and thus \cref{cor:b0_Mil1964_selfdual} yields 
\begin{equation*}
u(\svec(\calS_+^n)) =  (n^2 + n + 3) (2n^2 + 2n + 5)^{\frac{n^2 + n - 2}{2}}.
\end{equation*}
We expect this bound to be loose in view of \cref{ex:nno,ex:Lorentz}, but we emphasize that it is the first uniform upper bound for the cardinalities of $\svec(\calS_+^n)$-spectra of symmetric matrices.
This preliminary bound may serve as a starting point for deriving sharper uniform upper bounds.
\end{example}

\section{The exponential cone and its perturbed counterpart}\label{sec:nonsemialg}
As mentioned in \cref{sec:intro}, some cones are not symmetrically spectrally finite and hence are not definable in any o-minimal structure.
In this section, we provide an undefinable closed convex cone that is symmetrically spectrally finite but is not uniformly symmetrically spectrally finite.
Such a cone is constructed by perturbing the exponential cone, which is definable (see \cref{subsec:omin}) and thus uniformly symmetrically spectrally finite.
In \cref{subsec:expcone}, we provide an explicit uniform bound for the cardinalities of spectra of symmetric matrices relative to the exponential cone.
This illustrates the uniform symmetric spectral finiteness of definable cones shown in \cref{thm:omin_ussf}.
In \cref{subsec:nonsemialg_finite_nonunif}, using the discussion in \cref{subsec:expcone}, we provide an undefinable closed convex cone $K$ that is symmetrically spectrally finite but is not uniformly symmetrically spectrally finite.

Note that the situation considered in \cref{subsec:nonsemialg_finite_nonunif} differs from that in \cite[Remark~1]{See2021} in that the cone constructed in \cite{See2021} depends on the prescribed lower bound, whereas ours does not.
In \cite{See2021}, the lower bound is fixed \emph{before} the cone is constructed.
For any given $M$, a polyhedral cone, which is uniformly symmetrically spectrally finite~\cite[Proposition~6.1]{See1999}, was constructed such that the associated cone spectrum of some symmetric matrix had cardinality at least $M$.
In contrast, in \cref{subsec:nonsemialg_finite_nonunif}, the lower bound is chosen \emph{after} the cone is fixed.
We construct a cone such that for any given $M$, there exists a symmetric matrix whose cone spectrum has cardinality at least $M$.

In the subsequent subsections, we use perspective cones generated by the exponential function and its perturbed counterpart.
For a differentiable function $\omega\colon \bbR\to\bbR$, we define $f_{\omega}\colon\bbR\to \bbR$ as
\begin{equation}
f_{\omega}(t) \coloneqq \rme^t + \omega(t). \label{eq:f_eq_et_plus_omega}
\end{equation}
Furthermore, for $\bm{A} \in \calS^3$, we define
\begin{equation}
g_{\bm{A}}(\bm{x}) \coloneqq \begin{pmatrix}
1\\
0\\
x_3
\end{pmatrix}^\top\bm{A}\begin{pmatrix}
x_1\\
1\\
x_2
\end{pmatrix}(1 + x_1^2 + x_2^2) - \begin{pmatrix}
x_1\\
1\\
x_2
\end{pmatrix}^\top\bm{A}\begin{pmatrix}
x_1\\
1\\
x_2
\end{pmatrix}(x_1 + x_2x_3) \label{eq:def_pM}
\end{equation}
for $\bm{x}\in \bbR^3$, and define
\begin{align}
\bm{v}_{\omega}(t) & \coloneqq (t,1,f_{\omega}(t))^\top, \nonumber\\
\bm{n}_{\omega}(t) &\coloneqq \bm{v}_{\omega}'(t) \times \bm{v}_{\omega}(t) = (-f_{\omega}'(t),tf_{\omega}'(t)-f_{\omega}(t),1)^\top,\nonumber\\
V_{\omega}(t) &\coloneqq \bm{v}_{\omega}(t)^\top\bm{v}_{\omega}(t) = t^2 + 1 + f_{\omega}(t)^2, \nonumber\\
\lambda_{\bm{A},\omega}(t) &\coloneqq \frac{\bm{v}_{\omega}(t)^\top\bm{A}\bm{v}_{\omega}(t)}{V_{\omega}(t)},\nonumber\\
D_{\bm{A},\omega}(t) &\coloneqq g_{\bm{A}}(t,f_{\omega}(t),f_{\omega}'(t)) = \bm{v}_{\omega}'(t)^\top\bm{A}\bm{v}_{\omega}(t)V_{\omega}(t) - \bm{v}_{\omega}(t)^\top\bm{A}\bm{v}_{\omega}(t)\bm{v}_{\omega}'(t)^\top\bm{v}_{\omega}(t) \label{eq:def_DNq}
\end{align}
for $t\in \bbR$, where $\bm{v}_{\omega}'(t) \times \bm{v}_{\omega}(t)$ denotes the cross product of $\bm{v}_{\omega}'(t)$ and $\bm{v}_{\omega}(t)$.
Note that direct calculations lead to
\begin{equation}
\lambda_{\bm{A},\omega}'(t) = \frac{2D_{\bm{A},\omega}(t)}{V_{\omega}(t)^2}. \label{eq:lambdaNq_diff}
\end{equation}

First, we consider the situation where $f_{\omega}$ defined in \eqref{eq:f_eq_et_plus_omega} satisfies the following assumption.

\begin{assumption}\label{asm:perturbed_exp}
$f_{\omega}$ is a convex function of class $C^1$ and is equal to the exponential function outside the interval $[0,1]$.
In addition, $f_{\omega}'$ is nonnegative on $\bbR$.
\end{assumption}

Let $\psi \colon \bbR \to \bbR$ be such that $f_{\psi}$ satisfies Assumption~\ref{asm:perturbed_exp}.
We define $f_{\psi}^{\pi} \colon \bbR^2 \to (-\infty,+\infty]$ as
\begin{equation*}
f_{\psi}^{\pi}(\bm{x}) \coloneqq \begin{cases}
x_2f_{\psi}(x_1/x_2) & (x_2 > 0),\\
f_{\psi}^{\infty}(x_1) & (x_2 = 0),\\
+\infty & (x_2 < 0).
\end{cases}
\end{equation*}
Here, $f_{\psi}^{\infty}\colon \bbR \to (-\infty,+\infty]$ is the recession function of $f_{\psi}$ and is given by
\begin{equation}
f_{\psi}^{\infty}(t) = \lim_{h\to+\infty}\frac{f_{\psi}(s + ht) - f_{\psi}(s)}{h}, \label{eq:formula_finfty}
\end{equation}
where $s$ is an arbitrary element in $\bbR$~\cite[Theorem~8.5]{Roc1970}.
Then the \emph{perspective cone} generated by $f_{\psi}$ is defined as
\begin{equation}
K_{f_{\psi}} \coloneqq \epi f_{\psi}^{\pi}. \label{eq:epifpi}
\end{equation}
Since $f_{\psi}$ is a real-valued convex function by Assumption~\ref{asm:perturbed_exp}, $K_{f_{\psi}}$ is a closed convex cone and can be expressed as
\begin{equation}
K_{f_{\psi}} = \{\bm{x} \in \bbR^3 \mid x_2 > 0,\ (\tfrac{x_1}{x_2},\tfrac{x_3}{x_2})^\top \in \epi f_{\psi}\} \cup \{(x_1,0,x_3)^\top \mid (x_1,x_3)^\top \in \epi f_{\psi}^{\infty}\}; \label{eq:epifpi_explicit}
\end{equation}
see, for example, \cite[Equation~(2.8)]{WLP20XX} and the proof of \cite[Proposition~2.3.\Rnum{2}]{Com2018}.
The following lemma summarizes some properties of $K_{f_{\psi}}$.

\begin{lemma}\label{lem:property_K}
The following statements hold.
\begin{enumerate}[label=(\roman*), ref=\roman*]
\item $K_{f_{\psi}}$ can be represented as \label{enum:K_expression}
\begin{align}
K_{f_{\psi}} &= \{\bm{x}\in \bbR^3 \mid x_2 > 0,\ x_2f_{\psi}(x_1/x_2) \le x_3 \} \cup \{-a\bm{e}_1 + b\bm{e}_3 \mid a,b\ge 0\} \label{eq:K_expression_original}\\
&= \{a\bm{v}_{\psi}(t) + b\bm{e}_3 \mid a >0,\ b\ge 0,\ t\in \bbR\} \cup \{-a\bm{e}_1 + b\bm{e}_3 \mid a,b\ge 0\}. \label{eq:K_expression_cone}
\end{align}
\item For any $t\in \bbR$, we have $K_{f_{\psi}}^* \cap \{\bm{v}_{\psi}(t)\}^\perp = \{a\bm{n}_{\psi}(t)\mid a\ge 0\}$, where $\{\bm{v}_{\psi}(t)\}^\perp$ denotes the space of elements in $\bbR^3$ that are orthogonal to $\bm{v}_{\psi}(t)$. \label{enum:nt}
\item The set of nonzero elements in $K_{f_{\psi}}$ is the union of $S_1,\dots,S_5$ defined as \label{enum:nonzero_K}
\begin{align*}
S_1 &\coloneqq \setint K_{f_{\psi}},\\
S_2 &\coloneqq \{a\bm{v}_{\psi}(t) \mid a > 0,\ t\in \bbR\},\\
S_3 &\coloneqq \{-a\bm{e}_1 + b\bm{e}_3 \mid a,b > 0\},\\
S_4 &\coloneqq \{-a\bm{e}_1 \mid a > 0\},\\
S_5 &\coloneqq \{b\bm{e}_3 \mid b > 0\}.
\end{align*}
\end{enumerate}
\end{lemma}

\begin{proof}
First, we derive the expression in \eqref{eq:K_expression_original}.
The first set on the right-hand side of \eqref{eq:epifpi_explicit} is the same as its counterpart in \eqref{eq:K_expression_original}.
In addition, by Assumption~\ref{asm:perturbed_exp} and the formula for $f_{\psi}^{\infty}$ shown in \eqref{eq:formula_finfty}, we see that
\begin{equation*}
f_{\psi}^{\infty}(t) = \begin{cases}
0 & (t \le 0),\\
+\infty & (t > 0)
\end{cases}
\end{equation*}
and $\epi f_{\psi}^{\infty} = \{\bm{x}\in \bbR^2 \mid x_1 \le 0 \text{ and }x_2 \ge 0\}$.
Therefore, the second set on the right-hand side of \eqref{eq:epifpi_explicit} is the same as its counterpart in \eqref{eq:K_expression_original}.

Second, we derive the expression in \eqref{eq:K_expression_cone}.
It is sufficient to show that the first set on the right-hand side of \eqref{eq:K_expression_original} is the same as its counterpart in \eqref{eq:K_expression_cone}.
We assume that $\bm{x}$ belongs to the first set on the right-hand side of \eqref{eq:K_expression_original}.
We note that $x_2 > 0$.
Let $t \coloneqq x_1/x_2$ and $b \coloneqq x_3 - x_2f_{\psi}(t) \ge 0$.
Then we have
\begin{equation*}
\bm{x} = (x_2t,x_2,x_2f_{\psi}(t) + b)^\top = x_2\bm{v}_{\psi}(t) + b\bm{e}_3
\end{equation*}
so this vector belongs to the first set on the right-hand side of \eqref{eq:K_expression_cone}.
Conversely, we assume that $\bm{x}$ belongs to the first set on the right-hand side of \eqref{eq:K_expression_cone}.
Then there exist $a > 0$, $b \ge 0$, and $t\in \bbR$ such that
\begin{equation*}
\bm{x} = a\bm{v}_{\psi}(t) + b\bm{e}_3 = (at,a,af_{\psi}(t)+b)^\top.
\end{equation*}
Since $x_3 - x_2f_{\psi}(x_1/x_2) = b \ge 0$, we see that $\bm{x}$ belongs to the first set on the right-hand side of \eqref{eq:K_expression_original}.

Next, we prove \eqref{enum:nt}.
To prove that every nonnegative multiple of $\bm{n}_{\psi}(t)$ belongs to $K_{f_{\psi}}^* \cap \{\bm{v}_{\psi}(t)\}^\perp$, it is sufficient to show that $\bm{n}_{\psi}(t) \in K_{f_{\psi}}^* \cap \{\bm{v}_{\psi}(t)\}^\perp$.
By the definition of $\bm{n}_{\psi}$, it follows that $\bm{n}_{\psi}(t)\in \{\bm{v}_{\psi}(t)\}^\perp$.
Let $\bm{x}\in K_{f_{\psi}}$ be arbitrary.
First, we assume that $\bm{x}$ belongs to the first set on the right-hand side of \eqref{eq:K_expression_original}.
Let $s \coloneqq x_1/x_2$ for brevity.
Then it follows from the choice of $\bm{x}$ and the convexity of $f_{\psi}$ specified in Assumption~\ref{asm:perturbed_exp} that
\begin{equation*}
\bm{n}_{\psi}(t)^\top\bm{x} = -x_1f_{\psi}'(t) + x_2(tf_{\psi}'(t) - f_{\psi}(t)) + x_3 \ge x_2(f_{\psi}(s) - f_{\psi}(t) - f_{\psi}'(t)(s-t)) \ge 0.
\end{equation*}
Second, we assume that $\bm{x}$ belongs to the second set on the right-hand side of \eqref{eq:K_expression_original}, i.e., $\bm{x} = -a\bm{e}_1 + b\bm{e}_3$ with $a,b \ge 0$.
Then we have
\begin{equation*}
\bm{n}_{\psi}(t)^\top\bm{x} = af_{\psi}'(t) + b \ge 0,
\end{equation*}
where we use $f_{\psi}'(t) \ge 0$ from Assumption~\ref{asm:perturbed_exp} to derive the inequality.
In either case, it follows that $\bm{n}_{\psi}(t)^\top\bm{x} \ge 0$.
Since $\bm{x}\in K_{f_{\psi}}$ is arbitrary, we have $\bm{n}_{\psi}(t) \in K_{f_{\psi}}^*$.
Therefore, we obtain $\bm{n}_{\psi}(t) \in K_{f_{\psi}}^* \cap \{\bm{v}_{\psi}(t)\}^\perp$.

Conversely, let $\bm{x}\in K_{f_{\psi}}^* \cap \{\bm{v}_{\psi}(t)\}^\perp$.
We define $h\colon \bbR \to \bbR$ as
\begin{equation*}
h(s) \coloneqq \bm{x}^\top\bm{v}_{\psi}(s) = x_1s + x_2 + x_3f_{\psi}(s).
\end{equation*}
For every $s\in \bbR$, we see that $h(s) \ge 0$ since $\bm{v}_{\psi}(s) \in K_{f_{\psi}}$ by \eqref{enum:K_expression} and $\bm{x}\in K_{f_{\psi}}^*$.
In addition, we have $h(t) = 0$ by the choice of $\bm{x}$.
Therefore, $h$ attains a global minimum at $s = t$, and we have $0 = h'(t) = x_1 + x_3f_{\psi}'(t)$, i.e., $x_1 = -x_3f_{\psi}'(t)$.
Substituting this into $h(t) = 0$, we have $x_2 = x_3(tf_{\psi}'(t) - f_{\psi}(t))$.
Since $\bm{e}_3 \in K_{f_{\psi}}$ by \eqref{enum:K_expression} and $\bm{x} \in K_{f_{\psi}}^*$, it follows that $0\le \bm{x}^\top\bm{e}_3 = x_3$.
Thus, we obtain $\bm{x} = x_3\bm{n}_{\psi}(t)$, which is a nonnegative multiple of $\bm{n}_{\psi}(t)$.

Next, we prove \eqref{enum:nonzero_K}.
Let $\bm{x}$ be a nonzero element in $K_{f_{\psi}}$.
First, we assume that $\bm{x}$ belongs to the first set on the right-hand side of \eqref{eq:K_expression_cone}.
Then there exist $a > 0$, $b\ge 0$, and $t\in \bbR$ such that
\begin{equation*}
\bm{x} = a\bm{v}_{\psi}(t) + b\bm{e}_3 = (at,a,af_{\psi}(t)+b)^\top.
\end{equation*}
If $b > 0$, then $\bm{x}$ belongs to $S_1$ since $\bm{x}$ strictly satisfies all the inequalities defining the first set on the right-hand side of \eqref{eq:K_expression_original}.
If $b = 0$, then we have $\bm{x} = a\bm{v}_{\psi}(t) \in S_2$.
Second, we assume that $\bm{x}$ belongs to the second set on the right-hand side of \eqref{eq:K_expression_cone}.
Then there exist $a,b\ge 0$ such that $\bm{x} = -a\bm{e}_1 + b\bm{e}_3$.
Note that at least one of $a$ and $b$ is positive since $\bm{x}$ is nonzero.
If $a,b > 0$, then we have $\bm{x} \in S_3$.
If $a > 0$ and $b = 0$, then we have $\bm{x}\in S_4$.
If $a = 0$ and $b > 0$, then we have $\bm{x} \in S_5$.
Conversely, it follows from \eqref{enum:K_expression} that every element in $\bigcup_{i=1}^5 S_i$ is a nonzero element in $K_{f_{\psi}}$.
\end{proof}

Recall that $\lambda$ is a $K_{f_{\psi}}$-eigenvalue of $\bm{A}\in \calS^3$ if there exists a nonzero $\bm{x}\in K_{f_{\psi}}$ such that
\begin{equation}
\bm{x} \perp (\bm{A}\bm{x} - \lambda\bm{x}) \in K_{f_{\psi}}^*. \label{eq:def_Kfpsi_eig}
\end{equation}
To provide an upper bound for the number of $K_{f_{\psi}}$-eigenvalues of symmetric matrices, for every $\bm{A}\in \calS^3$ and $i\in \{1,\dots,5\}$, we define
\begin{equation*}
\Sp_i(\bm{A},K_{f_{\psi}}) \coloneqq \{\lambda \in \bbR \mid \text{there exists } \bm{x}\in S_i \text{ satisfying \eqref{eq:def_Kfpsi_eig}}\},
\end{equation*}
where $S_i$ is defined in \eqref{enum:nonzero_K} of \cref{lem:property_K}.
Since $\Sp(\bm{A},K_{f_{\psi}}) = \bigcup_{i=1}^5 \Sp_i(\bm{A},K_{f_{\psi}})$ by \eqref{enum:nonzero_K} of \cref{lem:property_K}, we have
\begin{equation}
\#\Sp(\bm{A},K_{f_{\psi}}) \le \sum_{i=1}^5\#\Sp_i(\bm{A},K_{f_{\psi}}). \label{eq:abs_Sp_Kfpsi}
\end{equation}
The following lemma provides an upper bound on $\#\Sp_i(\bm{A},K_{f_{\psi}})$ for each $i\in \{1,\dots,5\}$.

\begin{lemma}\label{lem:Spi_Kfpsi}
Let $\bm{A} \in \calS^3$.
Then the following statements hold.
\begin{enumerate}[label=(\roman*), ref=\roman*]
\item $\#\Sp_1(\bm{A},K_{f_{\psi}}) \le 3$. \label{enum:Sp1}
\item $\#\Sp_2(\bm{A},K_{f_{\psi}}) \le \#D_{\bm{A},\psi}^{-1}(0)$. \label{enum:Sp2}
\item $\#\Sp_3(\bm{A},K_{f_{\psi}}) \le 2$. \label{enum:Sp3}
\item $\#\Sp_4(\bm{A},K_{f_{\psi}}) \le 1$. \label{enum:Sp4}
\item $\#\Sp_5(\bm{A},K_{f_{\psi}}) \le 1$. \label{enum:Sp5}
\end{enumerate}
\end{lemma}

\begin{proof}
\eqref{enum:Sp1} follows from the fact that $K_{f_{\psi}}$-eigenvalues of $\bm{A}$ that admit $K_{f_{\psi}}$-eigenvectors belonging to $\setint K_{f_{\psi}}$ are eigenvalues of $\bm{A}$ in the usual sense~\cite[Remark~1.1]{Kuc1979}.

Next, we show \eqref{enum:Sp2}.
Let $\lambda \in \Sp_2(\bm{A},K_{f_{\psi}})$.
Then there exists $\bm{x}\in S_2$ satisfying \eqref{eq:def_Kfpsi_eig}.
For brevity, we let $\bm{s} \coloneqq \bm{A}\bm{x} - \lambda\bm{x}\in K_{f_{\psi}}^*$.
Since $\bm{x}\in S_2$, there exist $a > 0$ and $t\in \bbR$ such that $\bm{x} = a\bm{v}_{\psi}(t)$.
Then
\begin{equation}
\bm{s} = a(\bm{A}\bm{v}_{\psi}(t) - \lambda\bm{v}_{\psi}(t)), \label{eq:sigma4_l}
\end{equation}
and \eqref{eq:def_Kfpsi_eig} yields $0 = \bm{x}^\top\bm{s} = a^2(\bm{v}_{\psi}(t)^\top\bm{A}\bm{v}_{\psi}(t) - \lambda V_{\psi}(t))$.
This implies that $\lambda = \lambda_{\bm{A},\psi}(t)$.
In addition, by \eqref{enum:nt} of \cref{lem:property_K}, $\bm{s}$ is a nonnegative multiple of $\bm{n}_{\psi}(t)$.
It follows from the definition of $\bm{n}_{\psi}$ and \eqref{eq:sigma4_l} that
\begin{equation*}
0 = \bm{v}_{\psi}'(t)^\top\bm{s} = \frac{aD_{\bm{A},\psi}(t)}{V_{\psi}(t)},
\end{equation*}
and thus $D_{\bm{A},\psi}(t) = 0$.
Therefore, we obtain
\begin{equation*}
\#\Sp_2(\bm{A},K_{f_{\psi}}) \le \#\lambda_{\bm{A},\psi}(D_{\bm{A},\psi}^{-1}(0)) \le \#D_{\bm{A},\psi}^{-1}(0).
\end{equation*}

Next, we show \eqref{enum:Sp3}.
Let $\lambda \in \Sp_3(\bm{A},K_{f_{\psi}})$.
Then there exists $\bm{x}\in S_3$ satisfying \eqref{eq:def_Kfpsi_eig}.
For brevity, we let $\bm{s} \coloneqq \bm{A}\bm{x} - \lambda\bm{x} \in K_{f_{\psi}}^*$.
Since $\bm{x}\in S_3$, there exist $a,b > 0$ such that $\bm{x} = -a\bm{e}_1 + b\bm{e}_3$.
Then it follows from \eqref{eq:def_Kfpsi_eig} that $0 = \bm{x}^\top\bm{s} = -as_1 + bs_3$.
In addition, $-s_1$ and $s_3$ are nonnegative since $-\bm{e}_1,\bm{e}_3\in K_{f_{\psi}}$ by \eqref{enum:K_expression} of \cref{lem:property_K} and $\bm{s}\in K_{f_{\psi}}^*$.
Therefore, $s_1$ and $s_3$, the first and third entries of $\bm{A}\bm{x}-\lambda\bm{x}$, are actually $0$, and thus $\lambda$ is an eigenvalue of the matrix $\begin{psmallmatrix}
A_{11} & A_{13}\\
A_{13} & A_{33}
\end{psmallmatrix}$.
Since $2\times 2$ matrices have at most $2$ eigenvalues, we obtain $\#\Sp_3(\bm{A},K_{f_{\psi}}) \le 2$.

Next, we show \eqref{enum:Sp4} and \eqref{enum:Sp5}.
Let $i \in \{4,5\}$ and $\lambda \in \Sp_i(\bm{A},K_{f_{\psi}})$.
Then there exists $\bm{x} \in S_i$ satisfying \eqref{eq:def_Kfpsi_eig}.
Since $\bm{x} \in S_4\cup S_5$, there exist $a > 0$ and $\bm{u} \in \{-\bm{e}_1,\bm{e}_3\}$ such that $\bm{x} = a\bm{u}$.
Then \eqref{eq:def_Kfpsi_eig} implies that $\lambda = \bm{u}^\top\bm{A}\bm{u} / \norm{\bm{u}}^2$.
Therefore, we obtain $\#\Sp_i(\bm{A},K_{f_{\psi}}) \le 1$.
\end{proof}

\cref{lem:Spi_Kfpsi} implies that $K_{f_{\psi}}$ is (resp. uniformly) symmetrically spectrally finite if $\#D_{\bm{A},\psi}^{-1}(0)$, the number of zeros of $D_{\bm{A},\psi}$, is finite (resp. uniformly bounded with respect to $\bm{A}$).

Second, we consider the situation where $\omega$ is a polynomial function.

\begin{lemma}\label{lem:DAq_not_0}
Let $q\colon \bbR\to \bbR$ be a polynomial function.
In addition, let $\bm{A} \in \calS^3$ and suppose that $\bm{A} \not\in \bbR\bm{I}$.
Then $D_{\bm{A},q}$ is not identically $0$. 
\end{lemma}

\begin{proof}
We show this by contraposition.
We assume that $D_{\bm{A},q}$ is identically $0$.
It follows from \eqref{eq:lambdaNq_diff} and the assumption that $\lambda_{\bm{A},q}'$ is identically $0$.
This implies that there exists $\lambda \in \bbR$ such that $\lambda_{\bm{A},q}$ is identically $\lambda$.
In other words, 
\begin{align}
0 &= \bm{v}_q(t)^\top\bm{A}\bm{v}_q(t) - \lambda V_q(t) \nonumber\\
&= (A_{33}-\lambda)\rme^{2t} + 2((A_{33}-\lambda)q(t) + A_{13}t + A_{23})\rme^t +(A_{33}-\lambda)q(t)^2 + 2(A_{13}t+A_{23})q(t) \nonumber\\
&\quad + (A_{11}-\lambda)t^2 + 2A_{12}t + A_{22}-\lambda 
\label{eq:lambdaAq_const}
\end{align}
for all $t\in \bbR$.
Dividing \eqref{eq:lambdaAq_const} by $\rme^{2t}$ and letting $t\to +\infty$, we have $A_{33} = \lambda$.
With $A_{33} = \lambda$, dividing \eqref{eq:lambdaAq_const} by $\rme^t$ and letting $t\to +\infty$, we have $\lim_{t\to+\infty}(A_{13}t + A_{23}) = 0$.
This yields $A_{13} = A_{23} = 0$.
With $A_{33} = \lambda$ and $A_{13} = A_{23} = 0$, \eqref{eq:lambdaAq_const} implies that $(A_{11}-\lambda)t^2 + 2A_{12}t + A_{22}-\lambda = 0$ for all $t\in \bbR$.
Comparing the coefficients yields $A_{11} = A_{22} = \lambda$ and $A_{12} = 0$.
Therefore, we obtain $\bm{A} = \lambda\bm{I} \in \bbR\bm{I}$.
\end{proof}

\subsection{Uniform symmetric spectral finiteness of the exponential cone}\label{subsec:expcone}
The function $f_0$, obtained from \eqref{eq:f_eq_et_plus_omega} by taking the perturbation $\omega$ to be zero, is the exponential function and satisfies Assumption~\ref{asm:perturbed_exp}.
The perspective cone generated by $f_0$ as in \eqref{eq:epifpi} is the exponential cone defined in \eqref{eq:def_Kexp}.
In this subsection, we show that for any $\bm{A} \in \calS^3$, the number of $\Kexp$-eigenvalues of $\bm{A}$ is at most $18$.

First, we assume that $\bm{A}\in \bbR\bm{I}$.
Then it follows from \cref{lem:property_spec} that $\Sp(\bm{A},\Kexp)$ is a singleton.
Second, we assume that $\bm{A}\not\in \bbR\bm{I}$.
To apply \cref{lem:Spi_Kfpsi}, we bound the number of zeros of $D_{\bm{A},0}$, which is defined in \eqref{eq:def_DNq} with $\omega$ set to zero.

\begin{lemma}\label{lem:num_zeros_DA0}
Suppose that $\bm{A}\not\in \bbR\bm{I}$.
Then we have $\#D_{\bm{A},0}^{-1}(0) \le 11$.
\end{lemma}

\begin{proof}
For $t\in \bbR$, we have
\begin{align}
D_{\bm{A},0}(t) &= \underbrace{-A_{12}t^2 + (A_{11}-A_{22})t + A_{12}}_{\eqqcolon p_0(t)} \nonumber\\
&\quad+ \underbrace{(A_{13}t^3 + (-A_{13}+A_{23})t^2 + (A_{13}-2A_{23})t + A_{13}+A_{23})}_{\eqqcolon p_1(t)}\rme^t \nonumber\\
&\quad + \underbrace{((-A_{11}+A_{33})t^2 + (A_{11}-2A_{12}-A_{33})t + A_{12}-A_{22}+A_{33})}_{\eqqcolon p_2(t)}\rme^{2t} \nonumber\\
&\quad + \underbrace{(-A_{13}t + A_{13}-A_{23})}_{\eqqcolon p_3(t)}\rme^{3t}. \label{eq:DA0t}
\end{align}
Since $\bm{A}\not\in \bbR\bm{I}$, at least one of the polynomials $p_0,\dots,p_3$ is nonzero.
Equivalently, $J\coloneqq \{j\in \{0,1,2,3\} \mid p_j \neq 0\}$ is nonempty.
Let $d_0 \coloneqq 2$, $d_1 \coloneqq 3$, $d_2 \coloneqq 2$, and $d_3 \coloneqq 1$.
Then it follows from \cite[Problem~75 on page~46]{PS1998} that
\begin{equation*}
\#D_{\bm{A},0}^{-1}(0) \le - 1 + \sum_{j\in J}(\deg p_j +1) \le - 1 + \sum_{j\in J}(d_j +1) \le - 1 + \sum_{j = 0}^3(d_j +1) = 11,
\end{equation*}
which is the desired result.
\end{proof}

By \eqref{eq:abs_Sp_Kfpsi} and \cref{lem:Spi_Kfpsi,lem:num_zeros_DA0}, we obtain $\#\Sp(\bm{A},\Kexp) \le 18$.
We do not discuss here whether the upper bound for $\#\Sp(\bm{A},\Kexp)$ is tight because the main purpose of this subsection is to prove \cref{lem:num_zeros_DA0}, which is used in \cref{subsec:nonsemialg_finite_nonunif}.

\subsection{Example of an undefinable cone with finite but not uniformly finite spectra of symmetric matrices}\label{subsec:nonsemialg_finite_nonunif}
In this subsection, we provide a closed convex cone in $\bbR^3$ that is not definable in any o-minimal structure and is symmetrically spectrally finite but is not uniformly symmetrically spectrally finite.
We do this in the following three steps.
In \cref{subsubsec:construct_cone}, we construct a candidate for such a cone.
In \cref{subsubsec:finite_spectra}, we show that the cone is symmetrically spectrally finite.
In \cref{subsubsec:not_USSF}, we show that the cone is not uniformly symmetrically spectrally finite.

\subsubsection{Construction of a cone}\label{subsubsec:construct_cone}
For every positive integer $j$, we define $t_j \coloneqq 2^{-(j+1)}$, $w_j \coloneqq 2^{-(j+1)^2}$, $\alpha_j \coloneqq \frac{w_j}{4j}$, $l_j \coloneqq t_j - \frac{w_j}{2}$, $r_j \coloneqq t_j + \frac{w_j}{2}$, and $I_j \coloneqq [l_j,r_j]$.

\begin{lemma}\label{lem:Ij}
For every $j\ge 1$, we have
\begin{equation}
I_j \subseteq [\tfrac{7}{8}t_j,\tfrac{9}{8}t_j] \subseteq [0,1]. \label{eq:Ij_inclusion}
\end{equation}
In particular, for every $j,k \ge 1$ with $j + 1\le k$, we have $r_k < l_j$ and $I_j \cap I_k = \emptyset$.
\end{lemma}

\begin{proof}
For every $j\ge 1$, since $(j+1)^2 \ge j+3$, it follows that $w_j \le 2^{-(j+3)} = t_j/4$.
Since $0\le t_j \le 1/4$, we have \eqref{eq:Ij_inclusion}.
(As can be seen from this proof, $I_j$ is actually contained in $[0,9/32]$.
However, we use the larger interval $[0,1]$ for convenience in the subsequent analysis.)

In particular, for every $j,k \ge 1$ with $j + 1 \le k$, we see that
\begin{equation*}
l_j - r_k \ge \frac{7}{8}\cdot2^{-(j+1)} - \frac{9}{8}\cdot2^{-(k+1)} \ge 5\cdot 2^{-(j+5)} > 0,
\end{equation*}
i.e.,  $r_k < l_j$.
This implies that $I_j \cap I_k = \emptyset$.
\end{proof}

Let $\theta\colon \bbR \to \bbR$ be the function of period $4$ defined as
\begin{equation*}
\theta(t) \coloneqq \begin{cases}
t & (0\le t\le 1),\\
-t + 2 & (1\le t\le 3),\\
t - 4 & (3\le t\le 4).
\end{cases}
\end{equation*}
In addition, for every $j\ge 1$, we define $\phi_j\colon \bbR \to \bbR$ and $\psi_j\colon \bbR \to \bbR$ as
\begin{align*}
\phi_j(t) &\coloneqq \begin{cases}
\alpha_j\theta\left(\frac{t-l_j}{\alpha_j}\right) & (t\in I_j),\\
0 & (t\not\in I_j),
\end{cases}\\
\psi_j(t) &\coloneqq \int_{l_j}^t \phi_j(s)\,\rmd s.
\end{align*}
Figures~\ref{fig:phij} and \ref{fig:psij} show the graphs of $\phi_j$ and $\psi_j$, respectively.
The following lemma summarizes some properties of $\phi_j$ and $\psi_j$.

\begin{figure}[tbp]
  \centering
  \resizebox{0.95\linewidth}{!}{%
  \begin{tikzpicture}[
    x=0.5cm, y=0.5cm, >=latex,
    line cap=round, line join=round, font=\small,
    hguide/.style={black, densely dashed, thin},
    vguide/.style={black, densely dotted, thin}
  ]
    % 描画時には alpha_j = 2 とする
    \def\aval{2}

    % 省略区間の描画上の幅は 2alpha_j
    % 以下は省略後の描画用横座標
    \pgfmathsetmacro{\laststart}{6*\aval}
    \pgfmathsetmacro{\supportend}{10*\aval}
    \pgfmathsetmacro{\graphend}{11*\aval}
    \pgfmathsetmacro{\axisend}{12*\aval}

    % 横軸ラベルの上段の高さと上下段の間隔
    \pgfmathsetmacro{\labeltop}{-\aval-0.9}
    \def\labelstep{1.05}

    % y = alpha_j, -alpha_j の水平補助線
    % 実際の x = l_j + 4j alpha_j まで伸ばす
    \draw[hguide]
      (0,\aval) -- (\supportend,\aval);
    \draw[hguide]
      (0,{-\aval}) -- (\supportend,{-\aval});

    % 横軸ラベルと縦の補助線
    % q        ：描画用横座標を alpha_j で割った値
    % heightval：補助線の始点における関数値
    % labelrow ：0 は上段、1 は下段
    \foreach \q/\heightval/\labelrow/\ticktext in {
      0/0/0/{l_j},
      1/\aval/1/{l_j+\alpha_j},
      2/0/0/{l_j+2\alpha_j},
      3/-\aval/1/{l_j+3\alpha_j},
      4/0/0/{l_j+4\alpha_j},
      6/0/0/{l_j+4(j-1)\alpha_j},
      10/0/0/{l_j+4j\alpha_j}
    } {
      % ラベルの横方向の中央を指定した x 座標にそろえる
      \node[anchor=north, inner sep=2pt, outer sep=0pt]
        (xlabel\q)
        at ({\q*\aval},{\labeltop-\labelrow*\labelstep})
        {$\ticktext$};

      % ラベルの上辺中央まで縦の点線を引く
      \draw[vguide]
        ({\q*\aval},{\heightval}) -- (xlabel\q.north);

      % x軸上の目盛
      \draw[thin]
        ({\q*\aval},-0.12) -- ({\q*\aval},0.12);
    }

    % 最初と最後の1周期
    % 省略区間の内部には関数の線を描かない
    \foreach \s in {0,\laststart} {
      \draw[thick]
        (\s,0)
        -- ({\s+\aval},\aval)
        -- ({\s+3*\aval},{-\aval})
        -- ({\s+4*\aval},0);
    }

    % 台の外側では phi_j(x) = 0
    % 関数の描画終端は l_j + 4j alpha_j + alpha_j
    \draw[thick]
      ({-\aval},0) -- (0,0);
    \draw[thick]
      (\supportend,0) -- (\graphend,0);

    % x軸
    % 二重斜線との交点で分割し、その間には軸を描かない
    % 軸の右端は l_j + 4j alpha_j + 2 alpha_j
    \draw[thin, line cap=butt]
      ({-\aval},0) -- ({4.85*\aval},0);
    \draw[thin,->, line cap=butt]
      ({5.15*\aval},0) -- (\axisend,0)
      node[right] {$t$};

    % 省略を表す二重斜線
    \draw[thin]
      ({4.77*\aval},{-0.13*\aval})
      -- ({4.93*\aval},{0.13*\aval});
    \draw[thin]
      ({5.07*\aval},{-0.13*\aval})
      -- ({5.23*\aval},{0.13*\aval});

    % y軸：実際の x = l_j に配置
    \draw[thin,->]
      (0,{-\aval-0.45}) -- (0,{\aval+1});
      %node[above] {$y$};

    % y軸の目盛とラベル
    \draw[thin]
      (-0.12,\aval) -- (0.12,\aval);
    \node[left=3pt]
      at (0,\aval) {$\alpha_j$};

    \draw[thin]
      (-0.12,{-\aval}) -- (0.12,{-\aval});
    \node[left=3pt]
      at (0,{-\aval}) {$-\alpha_j$};

    % 関数名
    % \node[above]
    %   at ({8*\aval},{\aval+0.45})
    %   {$\phi_j(x)$};
  \end{tikzpicture}%
  }

\caption{Graph of $\phi_j$.
The function $\phi_j$ has period $4\alpha_j$ on $I_j = [l_j,l_j + 4j\alpha_j]$ in the sense that $\phi_j(t + 4k\alpha_j) = \phi_j(t)$ for all $t\in [l_j,l_j + 4\alpha_j]$ and $k \in \{0,\dots,j-1\}$, and it vanishes outside $I_j$.
}
  \label{fig:phij}
\end{figure}

\begin{figure}[tbp]
  \centering
  \resizebox{0.95\linewidth}{!}{%
  \begin{tikzpicture}[
    x=0.5cm, y=0.5cm, >=latex,
    line cap=round, line join=round, font=\small,
    hguide/.style={black, densely dashed, thin},
    vguide/.style={black, densely dotted, thin}
  ]
    % 描画時には alpha_j = 2 とする
    \def\aval{2}
    \pgfmathsetmacro{\asq}{\aval*\aval}
    \pgfmathsetmacro{\halfheight}{\asq/2}

    % 省略区間の描画上の幅は 2alpha_j
    % 以下は省略後の描画用横座標
    \pgfmathsetmacro{\laststart}{6*\aval}
    \pgfmathsetmacro{\supportend}{10*\aval}
    \pgfmathsetmacro{\graphend}{11*\aval}
    \pgfmathsetmacro{\axisend}{12*\aval}

    % 横軸ラベルの上段の高さと上下段の間隔
    \def\labeltop{-0.9}
    \def\labelstep{1.05}

    % y = alpha_j^2/2, alpha_j^2 の水平補助線
    % 実際の x = l_j + 4j alpha_j まで伸ばす
    \draw[hguide]
      (0,\halfheight) -- (\supportend,\halfheight);
    \draw[hguide]
      (0,\asq) -- (\supportend,\asq);

    % 横軸ラベルと縦の補助線
    % q        ：描画用横座標を alpha_j で割った値
    % heightval：補助線の始点における関数値
    % labelrow ：0 は上段、1 は下段
    \foreach \q/\heightval/\labelrow/\ticktext in {
      0/0/0/{l_j},
      1/\halfheight/1/{l_j+\alpha_j},
      2/\asq/0/{l_j+2\alpha_j},
      3/\halfheight/1/{l_j+3\alpha_j},
      4/0/0/{l_j+4\alpha_j},
      6/0/0/{l_j+4(j-1)\alpha_j},
      10/0/0/{l_j+4j\alpha_j}
    } {
      % ラベルの横方向の中央を指定した x 座標にそろえる
      \node[anchor=north, inner sep=2pt, outer sep=0pt]
        (xlabel\q)
        at ({\q*\aval},{\labeltop-\labelrow*\labelstep})
        {$\ticktext$};

      % ラベルの上辺中央まで縦の点線を引く
      \draw[vguide]
        ({\q*\aval},{\heightval}) -- (xlabel\q.north);

      % x軸上の目盛
      \draw[thin]
        ({\q*\aval},-0.12) -- ({\q*\aval},0.12);
    }

    % 最初と最後の1周期
    % 3つの二次関数を正確に表すベジェ曲線を使用する
    % 省略区間の内部には関数の線を描かない
    \foreach \s in {0,\laststart} {
      \begin{scope}[shift={(\s,0)}]
        \draw[thick]
          (0,0)

          % 0 <= u <= alpha_j：
          % psi_j = u^2/2
          .. controls ({\aval/3},0)
                  and ({2*\aval/3},{\asq/6})
          .. (\aval,\halfheight)

          % alpha_j <= u <= 3 alpha_j：
          % psi_j = alpha_j^2 - (u-2 alpha_j)^2/2
          .. controls ({5*\aval/3},{7*\asq/6})
                  and ({7*\aval/3},{7*\asq/6})
          .. ({3*\aval},\halfheight)

          % 3 alpha_j <= u <= 4 alpha_j：
          % psi_j = (u-4 alpha_j)^2/2
          .. controls ({10*\aval/3},{\asq/6})
                  and ({11*\aval/3},0)
          .. ({4*\aval},0);
      \end{scope}
    }

    % 台の外側では psi_j(x) = 0
    % 関数の描画終端は l_j + 4j alpha_j + alpha_j
    \draw[thick]
      ({-\aval},0) -- (0,0);
    \draw[thick]
      (\supportend,0) -- (\graphend,0);

    % x軸
    % 二重斜線との交点で分割し、その間には軸を描かない
    % 軸の右端は l_j + 4j alpha_j + 2 alpha_j
    \draw[thin, line cap=butt]
      ({-\aval},0) -- ({4.85*\aval},0);
    \draw[thin,->, line cap=butt]
      ({5.15*\aval},0) -- (\axisend,0)
      node[right] {$t$};

    % 省略を表す二重斜線
    \draw[thin]
      ({4.77*\aval},{-0.13*\aval})
      -- ({4.93*\aval},{0.13*\aval});
    \draw[thin]
      ({5.07*\aval},{-0.13*\aval})
      -- ({5.23*\aval},{0.13*\aval});

    % y軸：実際の x = l_j に配置
    \draw[thin,->]
      (0,-0.45) -- (0,{\asq+1});
    %   node[above] {$y$};

    % y軸の目盛とラベル
    \draw[thin]
      (-0.12,\halfheight) -- (0.12,\halfheight);
    \node[left=3pt]
      at (0,\halfheight) {$\frac12\alpha_j^2$};

    \draw[thin]
      (-0.12,\asq) -- (0.12,\asq);
    \node[left=3pt]
      at (0,\asq) {$\alpha_j^2$};

    % 関数名
    % \node[above]
    %   at ({8*\aval},{\asq+0.25})
    %   {$y=\psi_j(x)$};
  \end{tikzpicture}%
  }

\caption{Graph of $\psi_j$.
The function $\psi_j$ has period $4\alpha_j$ on $I_j = [l_j,l_j + 4j\alpha_j]$ in the sense that $\psi_j(t + 4k\alpha_j) = \psi_j(t)$ for all $t\in [l_j,l_j + 4\alpha_j]$ and $k \in \{0,\dots,j-1\}$, and it vanishes outside $I_j$.
}
\label{fig:psij}
\end{figure}

\begin{lemma}\label{lem:property_phij_etaj}
For every $j \ge 1$, the following statements hold.
\begin{enumerate}[label=(\roman*), ref=\roman*]
\item We have $\sup_{t\in \bbR}\abs{\phi_j(t)} = \alpha_j$ and $\sup_{t\in \bbR}\psi_j(t) = \alpha_j^2$. \label{enum:max_phij_etaj}
\item We have $0\le \psi_j(t) \le 1$ for all $t\in \bbR$. \label{enum:etaj_range}
\end{enumerate}
\end{lemma}

\begin{proof}
The following proof is visually supported by Figures~\ref{fig:phij} and \ref{fig:psij}.
Since $\phi_j$ satisfies
\begin{equation*}
\phi_j(t) = \begin{cases}
t-l_j & (l_j \le t \le l_j + \alpha_j),\\
-(t - l_j) +  2\alpha_j & (l_j + \alpha_j \le t \le l_j + 3\alpha_j),\\
t-l_j -4\alpha_j & (l_j + 3\alpha_j \le t \le l_j + 4\alpha_j),
\end{cases}
\end{equation*}
it attains its maximum $\alpha_j$ at $t = l_j + \alpha_j$ and its minimum $-\alpha_j$ at $t = l_j + 3\alpha_j$ on $[l_j,l_j + 4\alpha_j]$.
In addition, $\phi_j$ has period $4\alpha_j$ on $I_j$ in the sense that $\phi_j(t + 4k\alpha_j) = \phi_j(t)$ for all $t\in [l_j,l_j + 4\alpha_j]$ and $k \in \{0,\dots,j-1\}$.
Similarly, since $\psi_j$ satisfies
\begin{equation}
\psi_j(t) = \begin{cases}
\frac{1}{2}(t-l_j)^2 & (l_j \le t \le l_j + \alpha_j),\\
-\frac{1}{2}(t-(l_j+2\alpha_j))^2 + \alpha_j^2 & (l_j + \alpha_j \le t \le l_j + 3\alpha_j),\\
\frac{1}{2}(t-(l_j + 4\alpha_j))^2 & (l_j + 3\alpha_j \le t \le l_j + 4\alpha_j),
\end{cases} \label{eq:etaj_formula}
\end{equation}
it attains its maximum $\alpha_j^2$ at $t = l_j + 2\alpha_j$ and its minimum $0$ at $t \in \{l_j,l_j + 4\alpha_j\}$ on $[l_j,l_j + 4\alpha_j]$.
In addition, $\psi_j$ has period $4\alpha_j$ on $I_j$ in the sense that $\psi_j(t + 4k\alpha_j) = \psi_j(t)$ for all $t\in [l_j,l_j + 4\alpha_j]$ and $k \in \{0,\dots,j-1\}$.
Since $\phi_j$ and $\psi_j$ vanish outside $I_j$, and $\alpha_j^2 \le 1$, we obtain the desired results.
\end{proof}

By the disjointness of $(I_j)_{j\ge 1}$ shown in \cref{lem:Ij}, \eqref{enum:max_phij_etaj} of \cref{lem:property_phij_etaj}, and $\lim_{j\to\infty}\alpha_j = 0$, the series $\sum_{j=1}^{\infty}\phi_j$ and $\sum_{j=1}^{\infty}\psi_j$ converge uniformly on $\bbR$; denote their sums by $\phi$ and $\psi$, respectively.
Using $\psi$, we define $f_{\psi}$ as in \eqref{eq:f_eq_et_plus_omega}.
The following lemmas summarize some properties of $\phi$, $\psi$, and $f_{\psi}$.
In particular, we observe that $f_{\psi}$ satisfies Assumption~\ref{asm:perturbed_exp}.

\begin{lemma}\label{lem:property_f}
The following statements hold.
\begin{enumerate}[label=(\roman*), ref=\roman*]
\item $1\le f_{\psi}(t) \le \rme+1$ for all $t\in \bigcup_{j=1}^{\infty} I_j$. \label{enum:eta_nonnegative}
\item $f_{\psi}$ is of class $C^1$ and satisfies $f_{\psi}'(t) = \rme^t + \phi(t)$ for every $t\in \bbR$. \label{enum:f_C1}
\item For any $d \ge 0$, we have $\phi(t) = \rmo(t^d)$ and $\psi(t) = \rmo(t^d)$ as $t\downarrow 0$. \label{enum:phieta_oxm}
\item We have\label{enum:lim_ftminus1}
\begin{equation*}
\lim_{j\to\infty}\sup_{t\in I_j}(f_{\psi}(t) - 1) = 0.
\end{equation*}
\end{enumerate}
\end{lemma}

\begin{proof}
Let $t\in \bigcup_{j=1}^{\infty} I_j$.
Since $0 \le t \le 1$ by \cref{lem:Ij}, we have $1 \le \rme^t \le \rme$.
In addition, it follows from \eqref{enum:etaj_range} of \cref{lem:property_phij_etaj} that $0\le \psi(t) \le 1$.
Therefore, we obtain \eqref{enum:eta_nonnegative}.

For every $j\ge 1$, since $\phi_j$ is continuous, it follows from the fundamental theorem of calculus that $\psi_j$ is of class $C^1$ and satisfies $\psi_j' = \phi_j$.
By the uniform convergence of the series defining $\phi$ and $\psi$, it follows that $\psi$ is of class $C^1$ and satisfies $\psi' = \phi$.
Therefore, we obtain \eqref{enum:f_C1}.

Next, we prove that $\phi(t) = \rmo(t^d)$ as $t\downarrow 0$, as stated in \eqref{enum:phieta_oxm}.
For any $j \ge 1$ and $t\in I_j$, it follows from \cref{lem:Ij} and \eqref{enum:max_phij_etaj} of \cref{lem:property_phij_etaj} that
\begin{equation*}
\left|\frac{\phi(t)}{t^d}\right| = \left|\frac{\phi_j(t)}{t^d}\right| \le \frac{1}{4}\left(\frac{8}{7}\right)^d\frac{ 2^{-(j+1)(j+1-d)}}{j} \eqqcolon \bar{\phi}_j.
\end{equation*}
The sequence $(\bar{\phi}_j)_j$ is nonincreasing for $j \ge \frac{d}{2} - 1$ and satisfies $\lim_{j\to\infty} \bar{\phi}_j = 0$.
To show that $\lim_{t\downarrow 0} \phi(t)/t^d = 0$, take an arbitrary $\epsilon > 0$.
Let $j_0 \ge \frac{d}{2} - 1$ be a positive integer such that $\bar{\phi}_{j_0} < \epsilon$, and let $\delta \coloneqq r_{j_0}$.
Let $t \in (0,\delta)$.
If $t\not\in \bigcup_{j=1}^{\infty} I_j$, then we have $\abs{\phi(t)/t^d} = 0$.
If $t\in I_k$ for some $k\ge 1$, then since $k \ge j_0$ by \cref{lem:Ij}, we have
\begin{equation*}
\left|\frac{\phi(t)}{t^d}\right| \le \bar{\phi}_k \le \bar{\phi}_{j_0} < \epsilon.
\end{equation*}
In either case, it follows that $\abs{\phi(t)/t^d} < \epsilon$, and thus we obtain $\lim_{t\downarrow 0} \phi(t)/t^d = 0$.

Next, we prove that $\psi(t) = \rmo(t^d)$ as $t\downarrow 0$, as stated in \eqref{enum:phieta_oxm}.
For any $j \ge 1$ and $t\in I_j$, it follows from \cref{lem:Ij} and \eqref{enum:max_phij_etaj} of \cref{lem:property_phij_etaj} that
\begin{equation*}
0 \le \frac{\psi(t)}{t^d} = \frac{\psi_j(t)}{t^d} \le \frac{1}{16}\left(\frac{8}{7}\right)^d\frac{ 2^{-(j+1)(2j+2-d)}}{j^2} \eqqcolon \bar{\psi}_j.
\end{equation*}
Since the sequence $(\bar{\psi}_j)_j$ is nonincreasing for $j \ge \frac{d}{4} -1$ and $\lim_{j\to \infty}\bar{\psi}_j = 0$, in the same way as in the proof of $\lim_{t\downarrow 0} \phi(t)/t^d = 0$, we obtain $\lim_{t\downarrow 0} \psi(t)/t^d = 0$.

Next, \eqref{enum:lim_ftminus1} holds since
\begin{equation*}
0\le \sup_{t\in I_j}(f_{\psi}(t) - 1) \le \sup_{t\in I_j}(\rme^t-1) + \alpha_j^2 \le \rme^{\frac{9}{8}t_j} -1 + \alpha_j^2 \to 0\quad(j\to \infty),
\end{equation*}
where the first inequality follows from \eqref{enum:eta_nonnegative}, the second inequality holds by \eqref{enum:max_phij_etaj} of \cref{lem:property_phij_etaj}, the third inequality follows from \cref{lem:Ij} and the fact that the exponential function is nondecreasing, and the limit is derived from $\lim_{j\to\infty}\rme^{\frac{9}{8}t_j} = 1$ and $\lim_{j\to\infty}\alpha_j = 0$.
\end{proof}

\begin{lemma}
$f_{\psi}$ satisfies Assumption~\ref{asm:perturbed_exp}.
\end{lemma}

\begin{proof}
To begin with, we show that $f_{\psi}'$ is nondecreasing.
Let $t_1,t_2\in \bbR$ satisfy $t_1 \le t_2$.
First, we assume that $t_1\ge 0$ (and $t_2 \ge 0$).
It follows that
\begin{equation}
\rme^{t_2} - \rme^{t_1} \ge t_2 - t_1 \label{eq:exp_diff_ineq}
\end{equation}
since $\rme^t - t$ is nondecreasing for $t\ge 0$.
In addition, $\phi$ is $1$-Lipschitz continuous since the summands have pairwise disjoint supports by \cref{lem:Ij} and every $\phi_j$ for $j\ge 1$ is a piecewise linear function whose line segments and rays have slopes of $-1$, $0$, or $1$.
Therefore, we have
\begin{equation}
f_{\psi}'(t_2) - f_{\psi}'(t_1) = \rme^{t_2} - \rme^{t_1} + \phi(t_2) - \phi(t_1) \ge 0, \label{eq:fprime_x1ge0}
\end{equation}
where we use \eqref{enum:f_C1} of \cref{lem:property_f} to derive the equality and use \eqref{eq:exp_diff_ineq} and the $1$-Lipschitz continuity of $\phi$ to derive the inequality.
Second, we assume that $t_2 \le 0$ (and $t_1 \le 0$).
Since $f_{\psi}'(t) = \rme^t$ for every nonpositive $t$ by \eqref{enum:f_C1} of \cref{lem:property_f}, it follows from the nondecreasingness of the exponential function that
\begin{equation}
f_{\psi}'(t_2) - f_{\psi}'(t_1) = \rme^{t_2} - \rme^{t_1} \ge 0. \label{eq:fprime_x2le0} 
\end{equation}
Third, we assume that $t_1 < 0 < t_2$.
Then using \eqref{eq:fprime_x1ge0} and \eqref{eq:fprime_x2le0}, we have $f_{\psi}'(t_1) \le f_{\psi}'(0) \le f_{\psi}'(t_2)$.
In all three cases, it follows that $f_{\psi}'(t_1) \le f_{\psi}'(t_2)$, and thus $f_{\psi}'$ is nondecreasing.

We can now verify that $f_{\psi}$ satisfies Assumption~\ref{asm:perturbed_exp}.
The convexity of $f_{\psi}$ follows from the fact that $f_{\psi}'$ is nondecreasing.
That $f_{\psi}$ is of class $C^1$ was shown in \eqref{enum:f_C1} of \cref{lem:property_f}.
Since $\psi$ vanishes outside $\bigcup_{j=1}^\infty I_j$ by definition and this set is contained in $[0,1]$ by \cref{lem:Ij}, we see that $f_{\psi}$ is equal to the exponential function outside $[0,1]$.
The nonnegativity of $f_{\psi}'$ follows from the fact that $f_{\psi}'$ is nondecreasing and that $\lim_{t\to -\infty} f_{\psi}'(t) = \lim_{t\to -\infty}\rme^t = 0$.
Therefore, we obtain the desired result.
\end{proof}

The cone $K_{f_{\psi}}$ defined in \eqref{eq:epifpi} is the desired cone.
We denote it by $K$ hereafter.

\subsubsection{Symmetric spectral finiteness}\label{subsubsec:finite_spectra}
In this subsection, we show that the cone $K$ constructed in \cref{subsubsec:construct_cone} is symmetrically spectrally finite.
Let $\bm{A}\in \calS^3$.
If $\bm{A}\in \bbR\bm{I}$, it follows from \cref{lem:property_spec} that $\Sp(\bm{A},K)$ is a singleton.
In what follows, we assume that $\bm{A}\not\in\bbR\bm{I}$.
By \cref{lem:Spi_Kfpsi}, to prove that $\Sp(\bm{A},K)$ is finite, it is sufficient to show that the number of zeros of $D_{\bm{A},\psi}$, which is defined in \eqref{eq:def_DNq} with $\omega = \psi$, is finite.
We prove this in two steps.
First, we show in \cref{lem:DA_no_zero} that for some sufficiently small $\delta > 0$, the function $D_{\bm{A},\psi}|_{(0,\delta)}$ has no zeros.
Second, we show in \cref{lem:DA_delta1_finite_zeros} that the number of zeros of $D_{\bm{A},\psi}|_{\bbR \setminus (0,\delta)}$ is finite.

\begin{lemma}\label{lem:DA_no_zero}
Suppose that $\bm{A}\not\in\bbR\bm{I}$.
Then there exists $\delta > 0$ such that $D_{\bm{A},\psi}(t) \neq 0$ for all $t\in (0,\delta)$.
In other words, $D_{\bm{A},\psi}|_{(0,\delta)}$ has no zeros.
\end{lemma}

\begin{proof}
Since $g_{\bm{A}}$ defined in \eqref{eq:def_pM} is a polynomial, $g_{\bm{A}}$ is Lipschitz continuous on some neighborhood of $\bm{x} = (0,1,1)^\top$, and we write $L$ for a Lipschitz constant.
Since the two vector-valued functions $(t,f_{\psi}(t),f_{\psi}'(t))^\top$ and $(t,\rme^t,\rme^t)^\top$ are continuous and equal $(0,1,1)^\top$ at $t = 0$ by the definition of $f_{\psi}$ and \eqref{enum:f_C1} of \cref{lem:property_f}, for all $t$ in a sufficiently small neighborhood of $0$, we have
\begin{align}
\abs{D_{\bm{A},\psi}(t) - D_{\bm{A},0}(t)} &= \abs{g_{\bm{A}}(t,f_{\psi}(t),f_{\psi}'(t)) - g_{\bm{A}}(t,\rme^t,\rme^t)} \nonumber\\
&\le L\norm{(t,f_{\psi}(t),f_{\psi}'(t))^\top - (t,\rme^t,\rme^t)^\top} \nonumber\\
&\le L(\psi(t) + \abs{\phi(t)}). \label{eq:DA_minus_barDA_ineq}
\end{align}

Since $D_{\bm{A},0}$, given explicitly in \eqref{eq:DA0t}, has a Taylor series at the origin with infinite radius of convergence, this series represents the function on $\bbR$, so that
\begin{equation*}
D_{\bm{A},0}(t) = \sum_{n=0}^{\infty}\frac{D_{\bm{A},0}^{(n)}(0)}{n!}t^n 
\end{equation*}
for all $t\in \bbR$, where $D_{\bm{A},0}^{(n)}(0)$ denotes the $n$th derivative of $D_{\bm{A},0}$ at $0$.
If $D_{\bm{A},0}^{(n)}(0) = 0$ for all $n\ge 0$, then $D_{\bm{A},0}$ is identically $0$, which contradicts \cref{lem:DAq_not_0} under the assumption that $\bm{A} \not\in \bbR\bm{I}$.
Let $k$ be the minimum nonnegative integer $n$ such that $D_{\bm{A},0}^{(n)}(0) \neq 0$ and define $a \coloneqq D_{\bm{A},0}^{(k)}(0)/k! \neq 0$, so that
\begin{equation}
D_{\bm{A},0}(t) = at^k + \rmo(t^k)\quad(t\downarrow 0). \label{eq:barDA_atk_asymp}
\end{equation}

Therefore, we have
\begin{equation*}
D_{\bm{A},\psi}(t) = D_{\bm{A},0}(t) + \rmo(t^k) = at^k + \rmo(t^k) \quad(t\downarrow 0),
\end{equation*}
where the first equality follows from \eqref{eq:DA_minus_barDA_ineq} and \eqref{enum:phieta_oxm} in \cref{lem:property_f}, and the second from \eqref{eq:barDA_atk_asymp}.
This implies that $\lim_{t\downarrow 0}D_{\bm{A},\psi}(t)/t^k = a$, so there exists $\delta > 0$ such that
\begin{equation*}
\left|\frac{D_{\bm{A},\psi}(t)}{t^k} - a\right| \le \frac{\abs{a}}{2}
\end{equation*}
for all $t\in (0,\delta)$.
Thus, for any $t\in (0,\delta)$, we have $\abs{D_{\bm{A},\psi}(t)} \ge \abs{a}t^k/2 > 0$, and hence $D_{\bm{A},\psi}(t) \neq 0$.
\end{proof}

With $\delta$ as in \cref{lem:DA_no_zero}, we define $J_{\delta} \coloneqq \{j\ge 1 \mid I_j \cap [\delta,1] \neq \emptyset\}$.
We note that $J_{\delta}$ is finite since $\delta > 0$.
For every $j\in J_{\delta}$ and $k\in [4j]$, we define $I_{jk} \coloneqq [l_j + (k-1)\alpha_j,l_j + k\alpha_j]$.
We note that $I_j = \bigcup_{k=1}^{4j} I_{jk}$ for every $j\in J_{\delta}$.
For every $j\in J_{\delta}$ and $k\in [4j]$, by \eqref{eq:etaj_formula}, there exists a quadratic function $q_{jk}\colon \bbR \to \bbR$ such that $\psi|_{I_{jk}} = q_{jk}|_{I_{jk}}$.
This implies that $D_{\bm{A},\psi}|_{I_{jk}} = D_{\bm{A},q_{jk}}|_{I_{jk}}$.

\begin{lemma}\label{lem:DA_delta1_finite_zeros}
Suppose that $\bm{A}\not\in\bbR\bm{I}$.
Then the number of zeros of $D_{\bm{A},\psi}|_{\bbR\setminus (0,\delta)}$ is finite.
\end{lemma}

\begin{proof}
First, we show that the number of zeros of $D_{\bm{A},\psi}$ on $\bigcup_{j\in J_{\delta}} I_j$ is finite.
It suffices to show that for any $j \in J_{\delta}$ and $k \in [4j]$, the number of zeros of $D_{\bm{A},q_{jk}}|_{I_{jk}}$ is finite.
To derive a contradiction, assume that $D_{\bm{A},q_{jk}}|_{I_{jk}}^{-1}(0)$ is infinite for some $j \in J_{\delta}$ and $k \in [4j]$.
Since $D_{\bm{A},q_{jk}}|_{I_{jk}}^{-1}(0)$ is compact, there exists a sequence $(z_i)_{i\ge 1}$ of distinct elements in $D_{\bm{A},q_{jk}}|_{I_{jk}}^{-1}(0)$ that converges to some $z^* \in D_{\bm{A},q_{jk}}|_{I_{jk}}^{-1}(0)$.
Since $D_{\bm{A},q_{jk}}$ is analytic, the identity theorem~\cite[Corollary~1.2.7]{KP2002} implies that $D_{\bm{A},q_{jk}}$ is identically $0$.
This contradicts \cref{lem:DAq_not_0} under the assumption that $\bm{A} \not\in \bbR\bm{I}$.
Thus, $D_{\bm{A},q_{jk}}|_{I_{jk}}^{-1}(0)$ is finite for any $j \in J_{\delta}$ and $k \in [4j]$.

Second, we show that the number of zeros of $D_{\bm{A},\psi}$ on $\bbR \setminus ((0,\delta) \cup \bigcup_{j\in J_{\delta}}I_j)$ is finite.
For any $t\in \bbR \setminus ((0,\delta) \cup \bigcup_{j\in J_{\delta}}I_j)$, since $\psi(t)  = 0$ and $\phi(t) = 0$, it follows from the definition of $f_{\psi}$ and \eqref{enum:f_C1} of \cref{lem:property_f} that $f_{\psi}(t) = \rme^t$ and $f_{\psi}'(t) = \rme^t$.
This implies that $D_{\bm{A},\psi} = D_{\bm{A},0}$ on $\bbR \setminus ((0,\delta) \cup \bigcup_{j\in J_{\delta}}I_j)$.
It follows from \cref{lem:num_zeros_DA0} that the number of zeros of $D_{\bm{A},0}$ on $\bbR$ is finite.
Thus, the number of zeros of $D_{\bm{A},\psi}$ on $\bbR \setminus ((0,\delta) \cup \bigcup_{j\in J_{\delta}}I_j)$ is also finite.

Since the set $\bbR\setminus (0,\delta)$ is contained in the union of $\bigcup_{j\in J_{\delta}} I_j$ and $\bbR \setminus ((0,\delta) \cup \bigcup_{j\in J_{\delta}}I_j)$, the number of zeros of $D_{\bm{A},\psi}|_{\bbR\setminus (0,\delta)}$ is also finite.
\end{proof}

By \cref{lem:DA_no_zero,lem:DA_delta1_finite_zeros}, we conclude that the number of zeros of $D_{\bm{A},\psi}$ is finite, and thus $\Sp(\bm{A},K)$ is finite.

\subsubsection{Failure of uniform symmetric spectral finiteness}\label{subsubsec:not_USSF}
Let
\begin{equation*}
\bm{B}\coloneqq \begin{pmatrix}
-3 & -5 &  -4\\
-5 & -18&   8\\
-4 &   8&   0
\end{pmatrix} \quad\text{and}\quad \bm{C}\coloneqq \begin{pmatrix}
0 & 1 & 0\\
1 & 0&  0\\
0 & 0&  0
\end{pmatrix}.
\end{equation*}
We note that
\begin{align}
D_{\bm{B},0}(t) &= -5 + 4\rme^t + 13\rme^{2t} -12\rme^{3t} + t(15-20\rme^t + 7\rme^{2t} + 4\rme^{3t})\nonumber\\
&\quad+ t^2(5 + 12\rme^t + 3\rme^{2t}) - 4t^3\rme^t \nonumber\\
&= \frac{1}{3}t^4 + \rmo(t^4)\quad(t\downarrow 0), \label{eq:barDBt_Taylor}\\
D_{\bm{C},0}(t) &= 1 + \rme^{2t} -2t\rme^{2t} -t^2 = 2 + \rmo(1)\quad(t\downarrow 0). \label{eq:barDCt_Taylor}
\end{align}
In addition, for every $j\ge 1$, we define
\begin{align*}
\gamma_j &\coloneqq -\frac{D_{\bm{B},0}(t_j)}{D_{\bm{C},0}(t_j)} + 3\alpha_j,\\
\bm{A}_j &\coloneqq \bm{B} + \gamma_j\bm{C}.
\end{align*}
It follows from \eqref{eq:barDCt_Taylor} and $0\le t_j \le 1/4$ that $D_{\bm{C},0}(t_j) > 0$, so $\gamma_j$ is well-defined.
The asymptotic behavior of $\gamma_j$ is shown in the following lemma.

\begin{lemma}\label{lem:gammaj_Otj4}
We have $\gamma_j = \rmO(t_j^4)$ as $j\to \infty$.
In particular, $\lim_{j\to \infty} \gamma_j = 0$.
\end{lemma}

\begin{proof}
It follows from \eqref{eq:barDBt_Taylor} and \eqref{eq:barDCt_Taylor} that
\begin{align*}
D_{\bm{B},0}(t_j) &= \frac{1}{3}t_j^4 + \rmo(t_j^4)\quad(j\to \infty),\\
D_{\bm{C},0}(t_j) &= 2 + \rmo(1)\quad(j\to \infty).
\end{align*}
These estimates imply that $D_{\bm{B},0}(t_j)/D_{\bm{C},0}(t_j) = \rmO(t_j^4)$ as $j\to \infty$.
In addition, since
\begin{equation*}
\frac{\alpha_j}{t_j^4} = \frac{2^{-(j+1)(j-3)}}{4j} \to 0\quad(j\to\infty),
\end{equation*}
we have $\alpha_j = \rmo(t_j^4)$ as $j\to\infty$.
Therefore, it follows that $\gamma_j = \rmO(t_j^4)$ as $j\to\infty$.
Finally, since $\lim_{j\to \infty}t_j = 0$, we obtain $\lim_{j\to \infty} \gamma_j = 0$.
\end{proof}

In the following lemma, we establish a uniform approximation of $D_{\bm{A}_j,\psi}(t)$ by $6\alpha_j + 18\phi_j(t)$ for $t\in I_j$ and quantify the rate at which the approximation error tends to zero as $j\to\infty$.

\begin{lemma}\label{lem:supRjt}
For every $j\ge 1$, we define $R_j\colon \bbR\to \bbR$ as
\begin{equation*}
R_j(t) \coloneqq D_{\bm{A}_j,\psi}(t) - (6\alpha_j + 18\phi_j(t)).
\end{equation*}
Then we have
\begin{equation*}
\sup_{t \in I_j}\abs{R_j(t)} = \rmo(\alpha_j)\quad(j\to \infty).
\end{equation*}
\end{lemma}

We defer the proof of \cref{lem:supRjt} to Appendix~\ref{apdx:proof_supRjt} because it is lengthy and technical.
For every $j\ge 1$ and $k\in [j]$, we define
\begin{align*}
t_{jk}^+ &\coloneqq l_j + 4(k-1)\alpha_j + \alpha_j \in I_j,\\
t_{jk}^- &\coloneqq l_j + 4(k-1)\alpha_j + 3\alpha_j \in I_j.
\end{align*}
Since $\phi(t_{jk}^+) = \alpha_j$, we have $R_j(t_{jk}^+) = D_{\bm{A}_j,\psi}(t_{jk}^+) -24\alpha_j$.
Then it follows from \cref{lem:supRjt} that
\begin{equation*}
0 = \lim_{j\to \infty}\sup_{k\in [j]}\frac{\abs{R_j(t_{jk}^+)}}{\alpha_j} = \lim_{j\to \infty}\sup_{k\in [j]}\left|\frac{D_{\bm{A}_j,\psi}(t_{jk}^+)}{\alpha_j} -24\right|.
\end{equation*}
Therefore, there exists $j^+\ge 1$ such that $D_{\bm{A}_j,\psi}(t_{jk}^+) > 0$ for all $j\ge j^+$ and $k\in [j]$.
Similarly, since $R_j(t_{jk}^-) = D_{\bm{A}_j,\psi}(t_{jk}^-) + 12\alpha_j$, there exists $j^- \ge 1$ such that $D_{\bm{A}_j,\psi}(t_{jk}^-) < 0$ for all $j\ge j^-$ and $k\in [j]$.
For simplicity, we let
\begin{equation*}
\bar{j} \coloneqq \max\{j^+,j^-,2\}.
\end{equation*}
Since $D_{\bm{A}_j,\psi}(t)$ is continuous with respect to $t$, for any $j\ge\bar{j}$ and $k\in [j]$, there exists $\tau_{jk} \in (t_{jk}^+,t_{jk}^-)$ such that
\begin{equation}
D_{\bm{A}_j,\psi}(\tau_{jk}) = 0. \label{eq:DAjtaujk}
\end{equation}

We outline the two steps of the proof.
In the first step, we show in \cref{lem:lambdaAjtaujk_Keig} that $\lambda_{\bm{A}_j,\psi}(\tau_{jk})$ is a $K$-eigenvalue of $\bm{A}_j$ for any sufficiently large $j$ and any $k \in [j]$.
In the second step, we show in \cref{lem:lambdaAjtaujk_distinct} that for any sufficiently large $j$, the $K$-eigenvalues $\lambda_{\bm{A}_j,\psi}(\tau_{jk})$ are strictly increasing in $k\in [j]$.

\begin{lemma}\label{lem:lambdaAjtaujk_Keig}
There exists $j_0 \ge \bar{j}$ such that $\lambda_{\bm{A}_j,\psi}(\tau_{jk})$ is a $K$-eigenvalue of $\bm{A}_j$ for any $j\ge j_0$ and $k\in [j]$.
\end{lemma}

\begin{proof}
For brevity, we let $\bm{s}_{jk} \coloneqq \bm{A}_j\bm{v}_{\psi}(\tau_{jk}) - \lambda_{\bm{A}_j,\psi}(\tau_{jk})\bm{v}_{\psi}(\tau_{jk})$ for every $j\ge \bar{j}$ and $k\in [j]$.
Since $\bm{v}_{\psi}(\tau_{jk})$ is nonzero by definition, it is sufficient to show that $K \ni \bm{v}_{\psi}(\tau_{jk}) \perp \bm{s}_{jk} \in K^*$ for any sufficiently large $j$ and any $k\in [j]$.
Let $j\ge \bar{j}$ and $k\in [j]$.
Then $\bm{v}_{\psi}(\tau_{jk})\in K$ by \eqref{enum:K_expression} of \cref{lem:property_K}, and $\bm{v}_{\psi}(\tau_{jk}) \perp \bm{s}_{jk}$ by the definition of $\lambda_{\bm{A}_j,\psi}$.

In what follows, we shall show that $\bm{s}_{jk} \in K^*$ for any sufficiently large $j$ and any $k\in [j]$.
For every $j\ge \bar{j}$ and $k\in [j]$, direct calculations and \eqref{eq:DAjtaujk} imply that
\begin{equation*}
\bm{v}_{\psi}'(\tau_{jk})^\top\bm{s}_{jk} = \frac{D_{\bm{A}_j,\psi}(\tau_{jk})}{V_{\psi}(\tau_{jk})} = 0.
\end{equation*}
Therefore, $\bm{s}_{jk}$ is orthogonal to $\bm{v}_{\psi}(\tau_{jk})$ and $\bm{v}_{\psi}'(\tau_{jk})$.
By the definition of $\bm{n}_{\psi}$, there exists $a_{jk} \in \bbR$ such that $\bm{s}_{jk} = a_{jk}\bm{n}_{\psi}(\tau_{jk})$.

We define $\bm{s}_j\colon \bbR \to \bbR^3$ as $\bm{s}_j(t) \coloneqq \bm{A}_j\bm{v}_{\psi}(t) - \lambda_{\bm{A}_j,\psi}(t)\bm{v}_{\psi}(t)$ for every $j\ge \bar{j}$.
Similarly, we define $\bm{s}\colon \bbR \to \bbR^3$ as $\bm{s}(t) \coloneqq \bm{B}\bm{v}_{\psi}(t) - \lambda_{\bm{B},\psi}(t)\bm{v}_{\psi}(t)$.
Since
\begin{equation}
\sup_{t\in [0,1]}\norm{\bm{s}_j(t) - \bm{s}(t)} = \abs{\gamma_j}\sup_{t\in [0,1]}\norm{\bm{C}\bm{v}_{\psi}(t) - \lambda_{\bm{C},\psi}(t)\bm{v}_{\psi}(t)} \to 0\quad(j\to\infty) \label{eq:sup_sjt_minus_st}
\end{equation}
by \cref{lem:gammaj_Otj4}, we see that $(\bm{s}_j)_{j \ge \bar{j}}$ converges uniformly to $\bm{s}$ on $[0,1]$.
Moreover, it follows from \cref{lem:Ij} that
\begin{equation}
0\le \sup_{k\in [j]}\tau_{jk} \le \frac{9}{8}t_j \to 0\quad(j\to \infty). \label{eq:sup_taujk}
\end{equation}
Therefore, it follows from the triangle inequality, \eqref{eq:sup_sjt_minus_st}, \eqref{eq:sup_taujk}, and the continuity of $\bm{s}$ that
\begin{equation}
\sup_{k\in [j]}\norm{\bm{s}_j(\tau_{jk}) - \bm{s}(0)}
\le \sup_{k\in [j]}\norm{\bm{s}_j(\tau_{jk}) - \bm{s}(\tau_{jk})} + \sup_{k\in [j]}\norm{\bm{s}(\tau_{jk}) - \bm{s}(0)} \to 0\quad(j\to \infty). \label{eq:wjtaujk_minus_w0}
\end{equation}
Since
\begin{align*}
\bm{s}_j(\tau_{jk}) &= \bm{s}_{jk} =  a_{jk}\bm{n}_{\psi}(\tau_{jk}) = (-a_{jk}f_{\psi}'(\tau_{jk}),a_{jk}(\tau_{jk}f_{\psi}'(\tau_{jk})-f_{\psi}(\tau_{jk})),a_{jk})^\top,\\
\bm{s}(0) &= (-9,-9,9)^\top,
\end{align*}
\eqref{eq:wjtaujk_minus_w0} implies that
\begin{equation*}
\lim_{j\to\infty}\sup_{k\in [j]}\abs{a_{jk} - 9} = 0.
\end{equation*}
Thus, there exists $j_0 \ge \bar{j}$ such that $a_{jk} \ge 0$ for all $j\ge j_0$ and $k\in [j]$, and by \eqref{enum:nt} of \cref{lem:property_K}, we obtain $\bm{s}_{jk} =  a_{jk}\bm{n}_{\psi}(\tau_{jk}) \in K^*$.
\end{proof}

\begin{lemma}\label{lem:sigmajk}
For every $j\ge \bar{j}$ and $k\in [j-1]$, we define
\begin{equation*}
\sigma_{jk} \coloneqq \tau_{j,k+1} - \tau_{jk} - 4\alpha_j.
\end{equation*}
Then we have
\begin{equation*}
\sup_{k\in [j-1]}\abs{\sigma_{jk}} = \rmo(\alpha_j)\quad(j\to \infty).
\end{equation*}
In particular, it follows that
\begin{equation*}
\sup_{k\in [j-1]}\abs{\tau_{j,k+1}-\tau_{jk}} = \rmO(\alpha_j)\quad(j\to \infty).
\end{equation*}
\end{lemma}

\begin{proof}
Since $\phi_j(t) =  -t + l_j + 4(k-1)\alpha_j + 2\alpha_j$ for $t\in [t_{jk}^+,t_{jk}^-]$, we can substitute $t = \tau_{jk} \in [t_{jk}^+,t_{jk}^-]$ into the formula for $R_j(t)$ in \cref{lem:supRjt} and use \eqref{eq:DAjtaujk} to obtain

\begin{equation*}
\tau_{jk} = l_j + 4(k-1)\alpha_j + \frac{7}{3}\alpha_j + \frac{1}{18}R_j(\tau_{jk}).
\end{equation*}
Similarly, it follows that
\begin{equation*}
\tau_{j,k+1} = l_j + 4k\alpha_j + \frac{7}{3}\alpha_j + \frac{1}{18}R_j(\tau_{j,k+1}).
\end{equation*}
These equations yield
\begin{equation*}
\tau_{j,k+1} - \tau_{jk} = 4\alpha_j + \frac{1}{18}(R_j(\tau_{j,k+1}) - R_j(\tau_{jk})),
\end{equation*}
and we have
\begin{equation*}
\sup_{k\in [j-1]}\abs{\sigma_{jk}} = \frac{1}{18}\sup_{k\in [j-1]}\abs{R_j(\tau_{j,k+1}) - R_j(\tau_{jk})} \le \frac{1}{9}\sup_{t\in I_j}\abs{R_j(t)} = \rmo(\alpha_j)\quad(j\to\infty),
\end{equation*}
where we use the triangle inequality and $\tau_{j,k+1},\tau_{jk}\in I_j$ to derive the inequality and use \cref{lem:supRjt} to derive the final equality.
Therefore, we obtain $\sup_{k\in [j-1]}\abs{\sigma_{jk}} = \rmo(\alpha_j)$.

In particular, since
\begin{equation*}
\frac{1}{\alpha_j}\sup_{k\in [j-1]}\abs{\tau_{j,k+1}-\tau_{jk}} \le 4 + \frac{1}{\alpha_j}\sup_{k\in [j-1]}\abs{\sigma_{jk}} \to 4\quad(j\to \infty),
\end{equation*}
we obtain $\sup_{k\in [j-1]}\abs{\tau_{j,k+1}-\tau_{jk}} = \rmO(\alpha_j)$.
\end{proof}

\begin{lemma}\label{lem:yt}
Let
\begin{equation*}
y(t) \coloneqq \frac{1}{V_{\psi}(t)^2} - \frac{1}{4}.
\end{equation*}
Then we have
\begin{equation*}
\lim_{j\to\infty}\sup_{t\in I_j}\abs{y(t)} = 0.
\end{equation*}
\end{lemma}

\begin{proof}
For any $j\ge 1$ and $t\in I_j$, it follows that
\begin{equation}
4 \le V_{\psi}(t)^2 \le \left(1 + \frac{81}{64}t_j^2 + (\rme^{\frac{9}{8}t_j} + \alpha_j^2)^2\right)^2, \label{eq:Vt_eval}
\end{equation} 
where we use \eqref{enum:eta_nonnegative} of \cref{lem:property_f} to derive the first inequality and use \cref{lem:Ij} and \eqref{enum:max_phij_etaj} of \cref{lem:property_phij_etaj} to derive the second inequality.
Therefore, we have
\begin{equation*}
\sup_{t\in I_j}\abs{y(t)} \le \frac{1}{4} - \frac{1}{\left(1 + \frac{81}{64}t_j^2 + (\rme^{\frac{9}{8}t_j} + \alpha_j^2)^2\right)^2} \to 0\quad(j\to\infty),
\end{equation*}
where we use \eqref{eq:Vt_eval} to derive the inequality and use $\lim_{j\to\infty}t_j = \lim_{j\to \infty}\alpha_j = 0$ to derive the limit.
Therefore, we obtain the desired result.
\end{proof}

\begin{lemma}\label{lem:lambdaAjtaujk_distinct}
We have
\begin{equation*}
\lim_{j\to\infty}\sup_{k\in [j-1]}\left|\frac{\lambda_{\bm{A}_j,\psi}(\tau_{j,k+1}) - \lambda_{\bm{A}_j,\psi}(\tau_{jk})}{\alpha_j^2} -12\right| = 0.
\end{equation*}
In particular, there exists $j_1 \ge \bar{j}$ such that $\lambda_{\bm{A}_j,\psi}(\tau_{j,k+1}) > \lambda_{\bm{A}_j,\psi}(\tau_{jk})$ for any $j\ge j_1$ and $k\in [j-1]$.
\end{lemma}

\begin{proof}
For every $j\ge \bar{j}$, we see that
\begin{align}
&\sup_{k\in [j-1]}\left|\frac{\lambda_{\bm{A}_j,\psi}(\tau_{j,k+1}) - \lambda_{\bm{A}_j,\psi}(\tau_{jk})}{\alpha_j^2} -12\right|\nonumber\\
={}& \sup_{k\in [j-1]}\left|\frac{2}{\alpha_j^2}\int_{\tau_{jk}}^{\tau_{j,k+1}}\frac{D_{\bm{A}_j,\psi}(t)}{V_{\psi}(t)^2}\,\rmd t -12\right|\nonumber\\
={}& \sup_{k\in [j-1]}\left|\frac{2}{\alpha_j^2}\int_{\tau_{jk}}^{\tau_{j,k+1}}
(\tfrac{1}{4} + y(t))(6\alpha_j + 18\phi_j(t) + R_j(t))\,\rmd t -12\right|\nonumber\\
={}& \sup_{k\in [j-1]}\,\Biggl|\frac{3}{\alpha_j}\int_{\tau_{jk}}^{\tau_{j,k+1}}\rmd t - 12 + \frac{9}{\alpha_j^2}\int_{\tau_{jk}}^{\tau_{j,k+1}}\phi_j(t)\,\rmd t + \frac{2}{\alpha_j^2}\int_{\tau_{jk}}^{\tau_{j,k+1}}(\tfrac{1}{4} + y(t))R_j(t)\,\rmd t\nonumber\\
&\quad + \frac{2}{\alpha_j^2}\int_{\tau_{jk}}^{\tau_{j,k+1}}y(t)(6\alpha_j + 18\phi_j(t))\,\rmd t\Biggr|\nonumber\\
\le{}& \frac{3}{\alpha_j}\sup_{k\in [j-1]}\abs{\sigma_{jk}} + \frac{9}{\alpha_j^2}\sup_{k\in [j-1]}\left|\int_{\tau_{jk}}^{\tau_{j,k+1}}\phi_j(t)\,\rmd t\right| + \frac{2}{\alpha_j^2}\sup_{k\in [j-1]}\left|\int_{\tau_{jk}}^{\tau_{j,k+1}}(\tfrac{1}{4} + y(t))R_j(t)\,\rmd t\right|\nonumber\\
&\quad + \frac{2}{\alpha_j^2}\sup_{k\in [j-1]}\left|\int_{\tau_{jk}}^{\tau_{j,k+1}}y(t)(6\alpha_j + 18\phi_j(t))\,\rmd t\right|,\label{eq:frac_lambdadiff_alphaj2}
\end{align}
where we use \eqref{eq:lambdaNq_diff} to derive the first equality, use the definitions of $R_j$ and $y$ in \cref{lem:supRjt,lem:yt}, respectively, to derive the second equality, and use the definition of $\sigma_{jk}$ in \cref{lem:sigmajk} and the triangle inequality to derive the final inequality.
It is sufficient to show that all four terms on the right-hand side of \eqref{eq:frac_lambdadiff_alphaj2} converge to $0$ as $j\to \infty$.

First, it follows from \cref{lem:sigmajk} that 
\begin{equation*}
\lim_{j\to \infty}\frac{3}{\alpha_j}\sup_{k\in [j-1]}\abs{\sigma_{jk}} = 0.
\end{equation*}

Second, we see that
\begin{equation*}
\frac{9}{\alpha_j^2}\sup_{k\in [j-1]}\left|\int_{\tau_{jk}}^{\tau_{j,k+1}}\phi_j(t)\,\rmd t\right| = \frac{9}{\alpha_j^2}\sup_{k\in [j-1]}\left|\int_{\tau_{jk}+4\alpha_j}^{\tau_{j,k+1}}\phi_j(t)\,\rmd t\right| \le \frac{9}{\alpha_j}\sup_{k\in [j-1]}\abs{\sigma_{jk}} \to 0\quad(j\to\infty),
\end{equation*}
where the equality holds since $\int_{\tau_{jk}}^{\tau_{jk}+4\alpha_j}\phi_j(t)\,\rmd t = 0$ for every $k \in [j-1]$, the inequality holds by \eqref{enum:max_phij_etaj} of \cref{lem:property_phij_etaj} and the definition of $\sigma_{jk}$, and the limit is derived from \cref{lem:sigmajk}.

Third, it follows from the triangle inequality and \cref{lem:supRjt,lem:sigmajk,lem:yt} that
\begin{align*}
\frac{2}{\alpha_j^2}\sup_{k\in [j-1]}\left|\int_{\tau_{jk}}^{\tau_{j,k+1}}(\tfrac{1}{4} + y(t))R_j(t)\,\rmd t\right| &\le \frac{2}{\alpha_j^2}\left(\frac{1}{4}+\sup_{t\in I_j}\abs{y(t)}\right)\sup_{t\in I_j}\abs{R_j(t)}\sup_{k\in [j-1]}\abs{\tau_{j,k+1}-\tau_{jk}}\\
&\to 0\quad(j\to \infty).
\end{align*}

Fourth, we have
\begin{align*}
\frac{2}{\alpha_j^2}\sup_{k\in [j-1]}\left|\int_{\tau_{jk}}^{\tau_{j,k+1}}y(t)(6\alpha_j + 18\phi_j(t))\,\rmd t\right| &\le \frac{48}{\alpha_j}\sup_{t\in I_j}\abs{y(t)}\sup_{k\in [j-1]}\abs{\tau_{j,k+1}-\tau_{jk}}\\
&\to 0\quad(j\to\infty),
\end{align*}
where we use the triangle inequality and \eqref{enum:max_phij_etaj} of \cref{lem:property_phij_etaj} to derive the inequality and use \cref{lem:sigmajk,lem:yt} to derive the limit.
\end{proof}

We are ready to show that $K$ is not uniformly symmetrically spectrally finite.
Let $M > 0$ be arbitrary.
Recall that $j_0$ and $j_1$ are chosen in \cref{lem:lambdaAjtaujk_Keig,lem:lambdaAjtaujk_distinct}, respectively, and let $j \ge \max\{j_0,j_1,M\}$.
Since $j\ge j_0$ and $j\ge j_1$, it follows from \cref{lem:lambdaAjtaujk_Keig,lem:lambdaAjtaujk_distinct} that $\lambda_{\bm{A}_j,\psi}(\tau_{j1}),\dots,\lambda_{\bm{A}_j,\psi}(\tau_{jj})$ are distinct $K$-eigenvalues of $\bm{A}_j$.
Therefore, it follows from $j\ge M$ that $\#\Sp(\bm{A}_j,K) \ge j \ge M$, and thus $K$ is not uniformly symmetrically spectrally finite.
Note that the undefinability of $K$ in any o-minimal structure follows from \cref{thm:omin_ussf}.

\vspace{0.5cm}
\noindent
{\bf Acknowledgments}
The author is supported by JSPS Grant-in-Aid for Research Activity Start-up JP25K23344.

\vspace{0.5cm}
\noindent
{\bf Declarations}
This work was supported in part by a grant of access to OpenAI models through the ChatGPT for Academic Researchers program.
The author used ChatGPT to explore and discuss possible proof ideas.
The author subsequently developed and independently verified all proofs presented in this paper and takes full responsibility for their correctness.

\vspace{0.5cm}
\noindent
{\bf Data availability}
No data were used for the research described in the article.

\begin{appendices}
\section{Proof of \cref{lem:supRjt}}\label{apdx:proof_supRjt}
For every $j\ge 1$, it follows that
\begin{align}
\sup_{t\in I_j}\abs{R_j(t)} &= \sup_{t\in I_j}\abs{D_{\bm{A}_j,\psi}(t) - (6\alpha_j + 18\phi_j(t))} \nonumber\\
&\le \sup_{t\in I_j}\abs{D_{\bm{A}_j,\psi}(t) - (g_{\bm{A}_j}(t,f_{\psi}(t),\rme^t) + 18\phi_j(t))} \nonumber\\
&\quad + \sup_{t\in I_j}\abs{g_{\bm{A}_j}(t,f_{\psi}(t),\rme^t) - g_{\bm{A}_j}(t,\rme^t,\rme^t)} + \sup_{t\in I_j}\abs{D_{\bm{A}_j,0}(t) - 6\alpha_j}, \label{eq:supRjt_ineq}
\end{align}
where we use the triangle inequality and $D_{\bm{A}_j,0}(t) = g_{\bm{A}_j}(t,\rme^t,\rme^t)$ to derive the inequality.
To prove the desired result, it suffices to show that each of the three terms that appear on the right-hand side of \eqref{eq:supRjt_ineq} is $\rmo(\alpha_j)$ as $j\to \infty$.
\subsection{Asymptotic behavior of the first term}
Since $g_{\bm{A}}(\bm{x})$ defined in \eqref{eq:def_pM} is affine in $x_3$, it can be represented as
\begin{equation*}
g_{\bm{A}}(\bm{x}) = g_{\bm{A}}(x_1,x_2,0) + x_3(g_{\bm{A}}(x_1,x_2,1) - g_{\bm{A}}(x_1,x_2,0)).
\end{equation*}
In addition, $h_{\bm{A}}(x_1,x_2) \coloneqq g_{\bm{A}}(x_1,x_2,1) - g_{\bm{A}}(x_1,x_2,0)$ is a polynomial in $x_1$ and $x_2$ and is linear in $\bm{A}$.
Since
\begin{align*}
&D_{\bm{A}_j,\psi}(t) - (g_{\bm{A}_j}(t,f_{\psi}(t),\rme^t) + 18\phi_j(t))\\
={}& g_{\bm{A}_j}(t,f_{\psi}(t),\rme^t+\phi_j(t)) - (g_{\bm{A}_j}(t,f_{\psi}(t),\rme^t) + 18\phi_j(t))\\
={}& \phi_j(t)(h_{\bm{A}_j}(t,f_{\psi}(t)) - 18)
\end{align*}
for every $t\in I_j$, the first term on the right-hand side of \eqref{eq:supRjt_ineq} can be bounded by
\begin{align*}
&\sup_{t\in I_j}\abs{D_{\bm{A}_j,\psi}(t) - (g_{\bm{A}_j}(t,f_{\psi}(t),\rme^t) + 18\phi_j(t))}\\
={}& \sup_{t\in I_j}\abs{\phi_j(t)(h_{\bm{A}_j}(t,f_{\psi}(t)) - 18)}\\
\overset{\scriptsize \text{(a)}}\le{}& \alpha_j\sup_{t\in I_j}\abs{h_{\bm{A}_j}(t,f_{\psi}(t)) - 18}\\
\overset{\scriptsize \text{(b)}}={}& \alpha_j\sup_{t\in I_j}\abs{h_{\bm{A}_j}(t,f_{\psi}(t)) - h_{\bm{B}}(0,1)}\\
\overset{\scriptsize \text{(c)}}\le{}& \alpha_j\left(\sup_{t\in I_j}\abs{h_{\bm{A}_j}(t,f_{\psi}(t)) - h_{\bm{B}}(t,f_{\psi}(t))} + \sup_{t\in I_j}\abs{h_{\bm{B}}(t,f_{\psi}(t)) - h_{\bm{B}}(0,1)}\right),
\end{align*}
where we use \eqref{enum:max_phij_etaj} of \cref{lem:property_phij_etaj} to derive (a), use $h_{\bm{B}}(0,1) = 18$ to derive (b), and use the triangle inequality to derive (c).
It therefore suffices to show that
\begin{align}
\lim_{j\to \infty}\sup_{t\in I_j}\abs{h_{\bm{A}_j}(t,f_{\psi}(t)) - h_{\bm{B}}(t,f_{\psi}(t))} &= 0, \label{eq:1st_1_lim0}\\
\lim_{j\to \infty}\sup_{t\in I_j}\abs{h_{\bm{B}}(t,f_{\psi}(t)) - h_{\bm{B}}(0,1)} &= 0. \label{eq:1st_2_lim0}
\end{align}
We note that
\begin{equation}
(t,f_{\psi}(t))^\top \in \Omega \coloneqq [0,1]\times [1,\rme+1] \text{ for all $t\in \bigcup_{j=1}^{\infty}I_j$} \label{eq:def_setB}
\end{equation}
by \cref{lem:Ij} and \eqref{enum:eta_nonnegative} of \cref{lem:property_f}.

First, \eqref{eq:1st_1_lim0} holds since
\begin{equation*}
\sup_{t\in I_j}\abs{h_{\bm{A}_j}(t,f_{\psi}(t)) - h_{\bm{B}}(t,f_{\psi}(t))} = \abs{\gamma_j}\sup_{t\in I_j}\abs{h_{\bm{C}}(t,f_{\psi}(t))} \le \abs{\gamma_j}\sup_{\bm{x}\in \Omega}\abs{h_{\bm{C}}(\bm{x})} \to 0\quad(j\to \infty),
\end{equation*}
where the equality follows from the fact that $h_{\bm{A}}$ is linear in $\bm{A}$, the inequality is derived from \eqref{eq:def_setB}, and the limit follows from \cref{lem:gammaj_Otj4}.

Second, we show \eqref{eq:1st_2_lim0}.
Let $\epsilon > 0$ be arbitrary.
Since the polynomial $h_{\bm{B}}$ is uniformly continuous on $\Omega$, there exists $\delta > 0$ such that
for any $\bm{x},\bm{y}\in \Omega$ satisfying $\norm{\bm{x}-\bm{y}} < \delta$, we have
\begin{equation}
\abs{h_{\bm{B}}(\bm{x}) - h_{\bm{B}}(\bm{y})} < \epsilon. \label{eq:qB_unif_cont}
\end{equation}
Since $\sup_{t\in I_j} t$ and $\sup_{t\in I_j}(f_{\psi}(t) -1)$ converge to $0$ as $j\to \infty$ by \cref{lem:Ij} and \eqref{enum:lim_ftminus1} of \cref{lem:property_f}, respectively, there exists $j_2 \ge 1$ such that for any $j\ge j_2$ and $t\in I_j$, we have
\begin{equation}
0 < t < \frac{\delta}{2}  \quad\text{and}\quad 0 < f_{\psi}(t) -1 < \frac{\delta}{2}. \label{eq:eval_deltaover2}
\end{equation}
For any $j\ge j_2$ and $t\in I_j$, since $\norm{(t,f_{\psi}(t))^\top - (0,1)^\top} < \delta$ by \eqref{eq:eval_deltaover2} and $(t,f_{\psi}(t))^\top,(0,1)^\top \in \Omega$ by \eqref{eq:def_setB}, it follows from \eqref{eq:qB_unif_cont} that $\abs{h_{\bm{B}}(t,f_{\psi}(t))  - h_{\bm{B}}(0,1)} < \epsilon$.
Since $t\in I_j$ is arbitrary, we obtain
\begin{equation*}
\sup_{t\in I_j}\abs{h_{\bm{B}}(t,f_{\psi}(t))  - h_{\bm{B}}(0,1)} < \epsilon.
\end{equation*}
Therefore, \eqref{eq:1st_2_lim0} follows.

\subsection{Asymptotic behavior of the second term}
Let $\norm{\cdot}_{\rm F}$ denote the Frobenius norm.
Since $(\bm{A}_j)_{j\ge 1}$ is a bounded sequence by \cref{lem:gammaj_Otj4}, there exists $U > 0$ such that $\norm{\bm{A}_j}_{\rm F} \le U$ for all $j\ge 1$.
For every $j\ge 1$ and $t\in I_j$, by the mean value theorem, there exists $x_{jt}$ in $[\rme^t,f_{\psi}(t)]$, which is contained in $[1,\rme+1]$ by \cref{lem:Ij} and \eqref{enum:eta_nonnegative} of \cref{lem:property_f}, such that
\begin{equation*}
g_{\bm{A}_j}(t,f_{\psi}(t),\rme^t) - g_{\bm{A}_j}(t,\rme^t,\rme^t) = \partial_{x_2}g_{\bm{A}_j}(t,x_{jt},\rme^t)\psi_j(t),
\end{equation*}
where $\partial_{x_2}g_{\bm{A}}$ denotes the partial derivative of $g_{\bm{A}}$ with respect to $x_2$.
Therefore, we have
\begin{align*}
\sup_{t\in I_j}\abs{g_{\bm{A}_j}(t,f_{\psi}(t),\rme^t) - g_{\bm{A}_j}(t,\rme^t,\rme^t)} &= \sup_{t\in I_j}\abs{\partial_{x_2}g_{\bm{A}_j}(t,x_{jt},\rme^t)\psi_j(t)}\\
&\le \alpha_j^2 \sup_{t\in I_j}\abs{\partial_{x_2}g_{\bm{A}_j}(t,x_{jt},\rme^t)}\\
&\le \alpha_j^2 \sup_{\substack{t\in [0,1],\\x_2\in [1,\rme+1],\\\norm{\bm{A}}_{\rm F} \le U}}\abs{\partial_{x_2}g_{\bm{A}}(t,x_2,\rme^t)}\\
&=  \rmo(\alpha_j)\quad(j\to\infty),
\end{align*}
where the first inequality follows from \eqref{enum:max_phij_etaj} of \cref{lem:property_phij_etaj}, the second inequality is derived from \cref{lem:Ij}, $x_{jt} \in [1,\rme+1]$, and the choice of $U$, and the last equality holds since $\lim_{j\to \infty}\alpha_j = 0$.

\subsection{Asymptotic behavior of the third term}
For every $j\ge 1$ and $t\in I_j$, the mean value theorem gives $\kappa_{jt} \in I_j$ in the closed interval with endpoints $t$ and $t_j$ such that
\begin{align}
\abs{D_{\bm{A}_j,0}(t) - D_{\bm{A}_j,0}(t_j)} &= \abs{D_{\bm{A}_j,0}'(\kappa_{jt})(t-t_j)} \nonumber\\
&\le \frac{w_j}{2}\abs{D_{\bm{A}_j,0}'(\kappa_{jt})}\nonumber\\
&\le 2j\alpha_j\abs{D_{\bm{B},0}'(\kappa_{jt})} + 2j\alpha_j\abs{\gamma_jD_{\bm{C},0}'(\kappa_{jt})},\label{eq:DAjt_minus_DAjtj}
\end{align}
where the first inequality follows from $t\in I_j$ and the second follows from the triangle inequality and the fact that $D_{\bm{A},0}$ is linear in $\bm{A}$.
By the triangle inequality and \eqref{eq:DAjt_minus_DAjtj}, the third term on the right-hand side of \eqref{eq:supRjt_ineq} can be bounded by
\begin{equation}
\sup_{t\in I_j}\abs{D_{\bm{A}_j,0}(t) - 6\alpha_j} \le 2j\alpha_j\sup_{t\in I_j}\abs{D_{\bm{B},0}'(\kappa_{jt})} + 2j\alpha_j\abs{\gamma_j}\sup_{t\in I_j}\abs{D_{\bm{C},0}'(\kappa_{jt})} + \abs{D_{\bm{A}_j,0}(t_j) - 6\alpha_j}. \label{eq:3rd_ineq}
\end{equation}
It therefore suffices to show that all three terms on the right-hand side of \eqref{eq:3rd_ineq} are $\rmo(\alpha_j)$ as $j\to \infty$.

First, we show that $2j\alpha_j\sup_{t\in I_j}\abs{D_{\bm{B},0}'(\kappa_{jt})} = \rmo(\alpha_j)$.
It follows from \eqref{eq:barDBt_Taylor} that $D_{\bm{B},0}'(0) = 0$.
Combining this with the mean value theorem, we see that for every $j\ge 1$ and $t\in I_j$, there exists $s_{jt} \in (0,\kappa_{jt})$ satisfying $D_{\bm{B},0}'(\kappa_{jt}) = D_{\bm{B},0}''(s_{jt})\kappa_{jt}$.
Therefore, we have
\begin{equation*}
2j\alpha_j\sup_{t\in I_j}\abs{D_{\bm{B},0}'(\kappa_{jt})} = 2j\alpha_j\sup_{t\in I_j}(\abs{D_{\bm{B},0}''(s_{jt})}\kappa_{jt}) \le \frac{9}{4}j\alpha_jt_j\sup_{s\in [0,1]}\abs{D_{\bm{B},0}''(s)} = \rmo(\alpha_j)\quad(j\to \infty),
\end{equation*}
where we use $\kappa_{jt} \in I_j$ and \cref{lem:Ij} to derive the inequality and use $\lim_{j\to \infty}jt_j = 0$ to derive the final equality.

Second, $2j\alpha_j\abs{\gamma_j}\sup_{t\in I_j}\abs{D_{\bm{C},0}'(\kappa_{jt})} = \rmo(\alpha_j)$ holds since 
\begin{equation*}
2j\alpha_j\abs{\gamma_j}\sup_{t\in I_j}\abs{D_{\bm{C},0}'(\kappa_{jt})} \le 2j\alpha_j\abs{\gamma_j}\sup_{s\in [0,1]}\abs{D_{\bm{C},0}'(s)} = \rmo(\alpha_j) \quad(j\to \infty),
\end{equation*}
where the inequality follows from \cref{lem:Ij} and the equality holds since $\lim_{j\to \infty}j\abs{\gamma_j} = 0$ by \cref{lem:gammaj_Otj4} and $\lim_{j\to \infty}jt_j^4 = 0$.

Third, $\abs{D_{\bm{A}_j,0}(t_j) - 6\alpha_j} = \rmo(\alpha_j)$ holds since
\begin{equation*}
D_{\bm{A}_j,0}(t_j) = D_{\bm{B},0}(t_j) + \gamma_jD_{\bm{C},0}(t_j) = 3\alpha_jD_{\bm{C},0}(t_j) = 6\alpha_j + \rmo(\alpha_j)\quad(j\to \infty),
\end{equation*}
where the first equality holds since $D_{\bm{A},0}$ is linear in $\bm{A}$, the second equality holds by the definition of $\gamma_j$, and the third equality is a consequence of \eqref{eq:barDCt_Taylor}.
\end{appendices}

%-------------------
\bibliographystyle{plainurl} % スタイルの指定
% https://mathlandscape.com/latex-bibstyles/：【LaTeX】BibTeXのスタイル76個一覧
\bibliography{2026_3_ref} %
%-------------------
\end{document}